

\input epsf.tex

\def\2{{1\over 2}}

\def\d{\delta}
\def\a{\alpha}
\def\b{\beta}
\def\g{\gamma}

\def\e{\epsilon}
\def\l{\lambda}
\def\o{\omega}
\def\D{\Delta}

\def\fun#1#2#3{#1\colon #2\rightarrow #3}

\def\abs#1{\vert #1 \vert}
\def\frac#1#2{{{#1} \over {#2}}}

\def\st{\;\colon\;}
\def\tends{\rightarrow}
\def\weak{\rightharpoonup}

\def\dr{ {\rm d} }

\def\R{{\bf R}}
\def\N{{\bf N}}
\def\Z{{\bf Z}}

\def\T{{\bf T}}

\def\thm#1{\vskip 1 pc\noindent{\bf Theorem #1.\quad}\sl}
\def\lem#1{\vskip 1 pc\noindent{\bf Lemma #1.\quad}\sl}
\def\prop#1{\vskip 1 pc\noindent{\bf Proposition #1.\quad}\sl}
\def\cor#1{\vskip 1 pc\noindent{\bf Corollary #1.\quad}\sl}

\def\proof{\rm\vskip 1 pc\noindent{\bf Proof.\quad}}
\def\fin{\par\hfill $\backslash\backslash\backslash$\vskip 1 pc}
\def\txt#1{\quad\hbox{#1}\quad}

\def\L{{\cal L}}

\def\dcal{{\cal D}}
\def\G{\Gamma}

\def\o{\omega}
\def\r{\rho}

\def\cin#1{\2\abs{{#1}}^2}

\def\cinh#1{{1\over{2h}}\abs{{#1}}^2}

\def\2{\frac{1}{2}}
\def\inn#1#2{{\langle #1 ,#2\rangle}}
\def\Mt{{\cal M}_1(\T^p)}

\def\part{{\partial_{x}}}
\def\div{{\rm div}}

\def\pprime{{{}^\prime{}^\prime}}

\def\Ec{{\cal E}}
\def\gam{\{ \g_{\r,\tau} \}_{a\le\r\le\tau \le b}}



\baselineskip= 17.2pt plus 0.6pt
\font\titlefont=cmr17
\centerline{\titlefont A cost on paths of measures}
\vskip 1 pc
\centerline{\titlefont which induces the Fokker-Planck equation}
\vskip 4pc
\font\titlefont=cmr12
\centerline{         \titlefont {Ugo Bessi}\footnote*{{\rm 
Dipartimento di Matematica, Universit\`a\ Roma Tre, Largo S. 
Leonardo Murialdo, 00146 Roma, Italy.}}   }{}\footnote{}{
{{\tt email:} {\tt bessi@matrm3.mat.uniroma3.it}Work partially supported by the PRIN2009 grant "Critical Point Theory and Perturbative Methods for Nonlinear Differential Equations}} 
\vskip 0.5 pc
 
\par
\vskip 2pc
\centerline{\bf Abstract}
In [12], J. Feng and T. Nguyen define a cost on curves of measures which is finite exactly on the curves which solve a Fokker-Planck equation with $L^2$ drift. In this paper, using ideas of D. Gomes and E. Valdinoci, we give a different construction of the cost of [12].

\vskip 2 pc
\centerline{\bf  Introduction}
\vskip 1 pc

Let $\T^p\colon=\frac{\R^p}{\Z^p}$ denote the $p$-dimensional torus and let $\Mt$ be the space of Borel probability measures on 
$\T^p$, with the 2-Wasserstein distance $d_2$ on it. There is a notion of absolute continuity for paths valued in metric spaces such as $\Mt$; it turns out (see for instance theorem 8.3.1 of [2]) that a path $\fun{\mu}{(a,b)}{(\Mt,d_2)}$ is absolutely continuous if and only if it is a weak solution of the continuity equation, i. e. iff there is a vector field 
$X\in L^2((a,b)\times\T^p,\L^1\otimes\mu_t)$ (we shall always denote by $\L^p$ the $p$-dimensional Lebesgue measure) such that
$$\int_a^b\dr t\int_{\T^p}
[
\partial_t\phi+\inn{X}{\nabla_x\phi}
] \dr\mu_t(x)=0\qquad
\forall\phi\in C^\infty_0((a,b)\times\T^p)  .  $$
There is another way of defining the cost of a curve of measures in $\Mt$, which doesn't use metric spaces: we explain the approach of [12] for the Fokker-Planck equation, which is similar.  Let $\fun{\mu}{(a,b)}{\Mt}$ be a curve of measures; the cost of [12] is the term on the left in the equality below; the term on the right is its formal definition. 
$$\int_a^b\cin{\dot\mu_t-\2\D\mu_t}_{-1,\mu_t}\dr t=$$
$$\sup\{
\int_a^b\dr t\int_{\T^p}(-\partial_t\phi-\D\phi)\dr\mu_t(x)\st
\phi\in C^\infty_c((a,b)\times\T^p)\txt{and}
\int_a^b\dr t\int_{\T^p}|\nabla\phi|^2\dr\mu_t(x)\le 1
\}   .   \eqno (1)$$
If the $\sup$ above is finite, the Riesz representation theorem  yields the existence of a vector field 
$X\in L^2((a,b)\times\T^p,\L^1\otimes\mu_t)$ such that the weak form of the Fokker-Planck equation holds:
$$\int_a^b\dr t\int_{\T^p}
(-\partial_t\phi-\2\D\phi)\dr\mu_t(x)=
\int_a^b\dr t\int_{\T^p}\inn{\nabla\phi}{X}\dr\mu_t(x)
\qquad\forall\phi\in C^\infty_c((a,b)\times\T^p)   .  $$
Naturally, the cost (1) would be of little use if one did not prove that it is coercive and lower semicontinuous for the uniform convergence of curves of measures. Lower semicontinuity is immediate since (1) defines the cost as a $\sup$ of continuous functions. As for coercivity, in [12] it is  deduced from the following formula: if $\mu_t$ has density $\r_t$, then
$$\int_a^b\cin{\dot\r_t-\2\D\r_t}_{-1,\r_t}\dr t=
\int_a^b
\cin{\dot\r_t}_{-1,\r_t}\dr t+
\int_a^b\dr t\int_{\T^p}\2\frac{|\nabla\r_t|^2}{\r_t}
  \dr x   +
\int_{\T^p}[
\r_b\log\r_b-\r_a\log\r_a
]   \dr x  . $$
This formula connects together the entropy and 
$\frac{\nabla\r_t}{\r_t}$; following [15], we could call this latter object the "osmotic velocity". Indeed, $\mu_t=\r_t\L^p$ solves both the transport equation and Fokker-Planck; 
$\frac{\nabla\r_t}{\r_t}$ is the difference of the two drifts, i . e. is the component of the velocity due to osmosis. Being an integration by parts, the formula above requires some delicate estimates on 
$\r_t$. 

In this paper, we give an alternative definition of the cost (1) using the approximation scheme introduced in [13]. Namely, if 
$\mu\in\Mt$, we define 
${\cal D}_\mu$ as the set of the Borel functions 
$\fun{\g}{\T^p\times\R^p}{[0,+\infty)}$ such that $\g(x,\cdot)$ is a probability density on $\R^p$ for $\mu$ a. e. $x$; we interpret 
$\g(x,v)$ as the probability of jumping from $x$ to $x+v$ on the torus. If our particles are distributed with law $\mu\in\Mt$, we define $\mu\ast\g$ as their distribution after one jump (see section 1 below for the precise formula). 

Let $\mu_1,\mu_2\in\T^p$; according to [13] the cost of diffusing from $\mu_1$ at time $0$ to $\mu_2$ at time $h$  is  
$$\Ec^h(\mu_1,\mu_2)=
\inf    \left\{
\int_{\T^p\times\R^p}\left[
\cinh{v}+\log\g(x,v)
\right]  \g(x,v)\dr\mu_1(x)\dr v  
\st\g\in{\cal D}_{\mu_1},\quad
\mu_1\ast\g=\mu_2
\right\} -$$
$$\log\left(\frac{1}{2\pi h}\right)^\frac{p}{2}  .    \eqno (2)$$
If $\mu_t$ is a curve of measures defined on the interval $(a,b)$, we define
$$\Ec_{(a,b)}(\mu_t)=
\liminf_{h\tends 0}\frac{1}{h}
\int_a^{b-h}\Ec^h(\mu_t,\mu_{t+h})\dr t  .  $$
Heuristically, we are trying to replicate the fact that 
$u\in W^{1,2}(a,b)$ iff
$$\liminf_{h\tends 0}\int_a^{b-h}\left\vert
\frac{u(t+h)-u(t)}{h}
\right\vert^2\dr t<+\infty   .   $$
The last step in the construction of the cost is dictated by the fact that we want lower semicontinuity for free; thus, we shall set
$${\cal C}_{(a,b)}(\mu_t)=
\inf\liminf_{n\tends+\infty}{\cal E}_{(a,b)}(\mu^n_t)$$
where the $\inf$ is taken over all the sequences $\{ \mu^n_t \}$ which converge to $\mu_t$ uniformly. The reader should compare this with [12], where semicontinuity follows in a much more natural way.

We want to prove the following two theorems. 

\thm{1} Let $\fun{\mu}{(a,b)}{\Mt}$ be a Borel curve of measures. Let us suppose that the cost ${\cal C}_{(a,b)}(\mu_t)$ is finite. Then there is 
$X\in L^2((a,b)\times\T^p,\L^1\otimes\mu_t)$ such that $\mu_t$ is a weak solution of 
$$\partial_t\mu_t-\2\D\mu_t+\div(X\mu_t)=0
\quad t\in(a,b)    
\eqno (FP)_{X}  $$
Moreover, 
$${\cal C}_{(a,b)}(\mu_t)=
\int_a^b\int_{\T^p}|X(t,x)|^2\dr\mu_t(x)   .  \eqno (3)$$
Conversely, if $\mu_t$ is a weak solution of $(FP)_{X}$ with 
$X\in L^2((a,b)\times\T^p,\L^1\otimes\mu_t)$, then 
$${\cal C}_{(a,b)}(\mu_t)\le
\int_a^b\dr t\int_{\T^p}|X(t,x)|^2\dr\mu_t(x)  .  $$

\rm

\thm{2} Let us endow $\Mt$ with the 2-Wasserstein distance, which we shall call $d_2$ in the following. Then the two points below hold.

\noindent 1) If $\fun{\mu_n}{(a,b)}{\Mt}$ is a sequence of Borel curves such that, for some $M>0$, 
$${\cal C}_{(a,b)}(\mu_n)\le M\qquad\forall n\ge 1$$
then $\mu_n$ is compact in $C((a,b),\Mt)$.

\noindent 2) If $\fun{\mu_n}{(a,b)}{\Mt}$ is a sequence of Borel curves converging uniformly to $\fun{\mu}{(a,b)}{\Mt}$, then 
$${\cal C}_{(a,b)}(\mu)\le\liminf_{n\tends+\infty}
{\cal C}_{(a,b)}(\mu_n)  .  $$

\rm
\vskip 1pc

The paper is organized as follows: in section 1, we shall prove that, if the $\inf$ in (2) is finite, then it is a minimum and the  minimizer $\g$ is unique. This will allow us to define for $\mu_t$ a mean forward velocity and correlation matrix on the interval 
$[t,t+h]$; in section 2 we shall prove, in a tedious but elementary way, that cost, mean forward velocity and correlation matrix are Borel functions of $t$. In section 3, following [15], we let 
$h\tends 0$ and define the instantaneous forward velocity and correlation matrix; we prove that, if a path has finite cost, then the forward velocity is in $L^2$ and the correlation matrix is the identity. In section 4, we prove that paths of finite cost satisfy the Fokker-Planck equation; the drift of Fokker-Planck turns out to be the forward velocity defined in section 3. We also prove the $\ge$ half of equality (3). In section 5, we prove the converse: namely, we consider a semigroup $P_{s,t}$ on $\Mt$ induced by a Fokker-Planck equation with a smooth drift and show that the cost of the path $\mu_t\colon= P_{a,t}\mu_0$ is finite. Together with an approximation procedure, this will yield the $\le$ half of (3).  Most of the proofs depend on a few elementary estimates on the Gaussian; we have relegated them to the appendix.

The problem of finding solutions of the Fokker-Planck equation with irregular drift has been studied intensively; we refer the reader to [7], [8] and the bibliography therein for different approaches; [4] treats a problem strongly related to this. We also mention [1] and [14], two papers which connect the single step of Otto's scheme to large deviation theory. 

\vskip 1pc

\noindent{\bf Acknowledgement.} I would like to thank the referees for the stimulating comments.

\vskip 2pc
\centerline{\bf \S 1}
\centerline{\bf Definition of the cost}
\vskip 1pc

We begin with some notation, most of which is standard. 

\noindent $\bullet$) We define $\Mt$ as the space of all Borel probability measures on $\T^p$.

\noindent $\bullet$) We denote by 
$\fun{\pi}{\R^p}{\T^p\colon=\frac{\R^p}{\Z^p}}$ the natural projection; if $x,y\in\T^p$, we set
$$|x-y|_{\T^p}=\min\{
|\tilde x-\tilde y|\st\pi(\tilde x)=x,\quad\pi(\tilde y)=y
\}    .    $$

\noindent $\bullet$ Let $\G$ be a Borel probability measure on 
$\T^p\times\T^p$ whose first and second marginals are $\nu_0$ and $\nu_1$ respectively; we shall say that $\G$ is a transfer plan between $\nu_0$ and $\nu_1$. 

\noindent $\bullet$) For $\l\ge 1$, we define the $\l$-Wasserstein distance $d_\l$ on $\Mt$ by
$$d_\l(\nu_0,\nu_1)^\l=\min_\G
\int_{\T^p\times\T^p}|x-y|_{\T^p}^\l\dr\G(x,y)   \eqno (1.1)$$
where the minimum is over all transfer plans between  $\nu_0$ and $\nu_1$. It is standard ([2], [16]) that the minimum is attained and that $d_\l$ induces on $\Mt$ the weak$\ast$ topology; in particular, $\Mt$ is a compact metric space. In this paper, we shall use only $d_1$ and $d_2$.

\noindent $\bullet$) We denote by $Den$ the set of the Borel probability densities on $\R^p$.

\noindent $\bullet$) Let $\mu\in\Mt$; we say that a Borel function 
$\fun{\g}{\T^p\times\R^p}{\R}$ belongs to $\dcal_{\mu}$ if
$\g(x,\cdot)\in Den$ for $\mu$ a. e. $x\in\T^p$. 

\noindent $\bullet$) If $\mu\in\Mt$ and $\g\in\dcal_\mu$, we define the measure $\mu\ast\g$ on $\T^p$ by
$$\int_{\T^p}f(z)\dr(\mu\ast\g)(z)=
\int_{\T^p\times\R^p}f(x+v)\g(x,v)\dr\mu(x)\dr v$$
for every $f\in C(\T^p)$. Equivalently, $\mu\ast\g$ is the push-forward of the measure $\mu\otimes\g(x,\cdot)\L^p$ on 
$\T^p\times\R^p$ by the action of $\R^p$ on $\T^p$
$$\fun{}{\T^p\times\R^p}{\T^p},\qquad
\fun{}{(x,v)}{(x+v)}   .   $$
Heuristically, $\g(x,v)$ is the probability that a particle in $x$ will jump to $x+v$; if initially the particles are distributed with law 
$\mu$, after the jump they have law $\mu\ast\g$. 

\noindent $\bullet$) If $\mu_1,\mu_2\in\Mt$, we define 
$\dcal_{\mu_1,\mu_2}$ as the set of all $\g\in\dcal_{\mu_1}$ such that $\mu_1\ast\g=\mu_2$. 

\vskip 1pc

Since $\dcal_{\mu_1,\mu_2}$ could be empty, we need two lemmas. 

\lem{1.1} $\dcal_{\mu_1,\mu_2}\not=\emptyset$ iff $\mu_2<<\L^p$.

\proof Let $\g\in\dcal_{\mu_1,\mu_2}$; setting 
$\tilde\g(x,v)=\g(x,-x+v)$ we easily check that
$$\G=\mu_1\otimes\pi_\sharp(\tilde\g(x,\cdot)\L^p)   \eqno (1.2)$$
is a transfer plan between $\mu_1$ and $\mu_2$; conversely, if 
$\G$ defined as in (1.2) is a transfer plan between $\mu_1$ and 
$\mu_2$, then $\g\in\dcal_{\mu_1,\mu_2}$. Thus, it suffices to show that there is a transfer plan of the form (1.2) iff 
$\mu_2<<\L^p$.

We begin with the if part. If $\mu_2=\a(x)\L^p$ for some 
$\a\in L^1(\T^p)$, we set 
$$\tilde\g(x,v)=\a(\pi(v))1_{[0,1)^p}(v)  .  $$
Now it is easy to see that $\G$ defined as in (1.2) is a transfer plan between $\mu_1$ and $\mu_2$.

Conversely, if $\G$ is a transfer plan from $\mu_1$ to $\mu_2$ and $f\in C(\T^p,\R)$, we have the first equality below; the second one holds if $\G$ has the form (1.2), while the third one follows by the periodicity of $f$. 
$$\int_{\T^p}f(y)\dr\mu_2(y)=
\int_{\T^p\times\T^p}f(y)\dr\G(x,y)=
\int_{\T^p\times\R^p}f(v)\tilde\g(x,v)\dr\mu_1(x)\dr v=$$
$$\sum_{k\in\Z^p}\int_{\T^p\times[0,1)^p}
f(v)\tilde\g(x,v+k)\dr\mu_1(x)\dr v  .  $$
By dominated convergence this means that $\mu_2$, considered as a measure on $[0,1)^p$, has density
$$\sum_{k\in\Z^p}\int_{\T^p}\tilde\g(x,v+k)\dr\mu_1(x)  $$
ending the proof.

\fin

\lem{1.2} 1) Let $\mu_1\in\Mt$. Then, the map
$$\fun{\Psi}{{\cal D}_{\mu_1}}{\Mt},\qquad
\fun{\Psi}{\g}{\mu_1\ast\g}  $$
is continuous from the weak topology of 
$L^1(\T^p\times\R^p,\mu_1\otimes\L^p)$ to the weak$\ast$ topology of $\Mt$.

\noindent 2) $\dcal_{\mu_1,\mu_2}$ is convex and weakly closed in $L^1(\mu_1\otimes\L^p)$.

\proof We prove point 1). Let $f\in C(\T^p)$; the formula below shows that $\Psi$ sends a weak neighbourhood of $\g$  into a weak$\ast$ neighbourhood of $\mu_1\ast\g$; the equality is  the definition of $\mu_1\ast\g$.
$$\left\{
\g^\prime\st
\left\vert
\int_{\T^p}f(x)\dr(\mu_1\ast\g)(x)-
\int_{\T^p}f(x)\dr(\mu_1\ast\g^\prime)(x)
\right\vert   <\e
\right\}   =  $$
$$\left\{
\g^\prime\st
\left\vert
\int_{\T^p\times\R^p}f(x+y)[\g(x,y)-\g^\prime(x,y)]
\dr\mu_1(x)\dr y
\right\vert   <\e
\right\}   .  $$

As for point 2), the fact that $\dcal_{\mu_1,\mu_2}$ is convex follows because the map $\fun{}{\g}{\mu\ast\g}$ is linear. As for the weak closure, we note that, by point 1), $\dcal_{\mu_1,\mu_2}$ is relatively closed in $\dcal_{\mu_1}$; thus, it suffices to show that 
$\dcal_{\mu_1}$ is weakly closed. By the definition of 
${\cal D}_{\mu_1}$, we have to to show that, if $\g$ is in the weak closure of $\dcal_{\mu_1}$, then $\g$ is positive (which is standard) and 
$$\int_{\R^p}\g(x,v)\dr v=1\txt{for $\mu_1$ a. e.} x\in\T^p  .  $$
In turn, this is implied by 
$$\mu_1(B)=\int_{\T^p\times\R^p}
\g(x,v)1_B(x)\dr\mu_1(x)\dr v \eqno (1.3)$$
for all Borel sets $B\subset\T^p$. Since the functional
$$\fun{}{\g}{
\int_{\T^p\times\R^p}1_B(x)\g(x,v)\dr\mu_1(x)\dr v
}  $$
is continuous for the weak topology of $L^1(\mu_1\times\L^p)$, this follows from the fact that (1.3) holds for the functions of 
$\dcal_{\mu_1}$. 

\fin

Our cost has been introduced in [13]: roughly, it is the kinetic energy minus the entropy: in dynamical terms, the pressure of the kinetic energy.

\vskip 1pc
\noindent{\bf Definition.} For $h>0$ and 
$\fun{\g}{\T^p\times\R^p}{[0,+\infty)}$, we define
$$A_h(\g,(x,v))=\cinh{v}\g(x,v)+\g(x,v)\log\g(x,v)  .  $$
If $\g$ depends only on $v$, we shall write $A_h(\g,v)$.

\vskip 1pc

We recall a few facts from [6] about the functional of Gomes and Valdinoci.

\lem{1.3} Let $h>0$ and let $\mu_1\in\Mt$. Then, the following points hold.

\noindent 1) The functional
$$\fun{I}{\dcal_{\mu_1}}{\R\cup\{ +\infty \}},\qquad
\fun{I}{\g}{\int_{\T^p\times\R^p}A_h(\g,(x,v))\dr\mu_1(x)\dr v}$$
is well defined.

\noindent 2) The functional $I$ is l. s. c. for the weak topology of 
$L^1(\mu_1\otimes\L^p)$.

\noindent 3) If $M\in\R$, the set
$$E_M=\{
\g\in\dcal_{\mu_1}\st I(\g)\le M
\}  $$
is uniformly integrable for the measure $\mu_1\otimes\L^p$.

\noindent 4) There is an increasing function 
$\fun{B}{\R}{[0,+\infty)}$, independent of $\mu_1$, such that, if 
$\g\in E_M$, then
$$\int_{\T^p\times\R^p}\cinh{v}\g(x,v)\dr\mu_1(x)\dr v\le
B(M) .  \eqno (1.4)$$

\noindent 5) The set $E_M$ is weakly compact in 
$L^1(\mu_1\otimes\L^p)$.

\noindent 6) Let $\mu_2\in\Mt$  be such that 
$\dcal_{\mu_1,\mu_2}$ is not empty. Then, the functional $I$ has a unique minimum 
$\g_{\mu_1,\mu_2}$ in $\dcal_{\mu_1,\mu_2}$.

\proof We only sketch the proof of this lemma and refer the reader to [6] for details. We begin with point 1). Note that
$$\cinh{v}y+y\log y\ge-e^{-1}e^{-\cinh{v}}\qquad
\forall y\ge 0$$
and thus
$$A_h(\g,(x,v))\ge-e^{-1}e^{-\cinh{v}}   .   \eqno (1.5)$$
Since the term on the right belongs to $L^1(\mu\otimes\L^p)$, we get that the integral of the negative part of $A_h(\g,(x,v))$ is finite; thus the integral of $A_h(\g,(x,v))$ is well defined, though possibly $+\infty$. 

We prove point 2). Since $I$ is convex, it suffices to prove that $I$ is l. s. c. for the strong topology of $L^1(\mu\otimes\L^p)$; as in lemma 1.3 of [6], this follows from Fatou's lemma and (1.5).

As for point 3), it follows as in lemma 1.2 of [6]: grossly, uniform integrability follows since $\g\log\g$ is superlinear. 

Again referring to lemma 1.2 of [6] for details, we derive (1.4) from lemma A.1 of the appendix in the following way. We set
$$a(x)=\frac{1}{h}\int_{\R^p}v\g(x,v)\dr v, \qquad
\d(x)=\frac{1}{h}\int_{\R^p}|v-ha(x)|^2\g(x,v)\dr v  .  $$
It is easy to see that point 4) follows if we prove that there is a function $C(M)$ such that 
$$\frac{h}{2}\int_{\T^p}|a(x)|^2\dr\mu_1(x)+
\int_{\T^p}\d(x)\dr\mu_1(x)\le C(M)  \eqno (1.6)$$
for all $\g\in E_M$. Let us set 
$$F=\{
x\in\T^p\st \d(x)\ge 2
\}  .  $$
The first inequality below comes since $\g\in E_M$, the second one comes from lemma A.1, the third one from formula (A.10) of the appendix and the definition of $T(a,\d)$.
$$M\ge\int_{\T^p\times\R^p}
A_h(\g,(x,v))\dr\mu_1(x)\dr v\ge
\int_{\T^p}T(a(x),\d(x))\dr\mu_1(x)=$$
$$\int_{\T^p}[T(a(x),\d(x))-T(a(x),1)]\dr\mu_1(x)+
\int_{\T^p}T(a(x),1)\dr\mu_1(x)    \ge$$
$$\int_FD_1(\d(x)-1)\dr\mu_1(x)+
\int_{F^c}D_1|\d(x)-1|^2\dr\mu_1(x)+
\int_{\T^p}\frac{h}{2}|a(x)|^2\dr\mu_1(x)+
\log\left(\frac{1}{2\pi h}\right)^\frac{p}{2}  .  $$
By H\"older's inequality, this implies (1.6). 

To prove point 5), we note that $E_M$ is weakly closed by point 2) of this lemma. Thus, it suffices to prove that it is relatively compact; this follows  by point 3) of this lemma and the fact that the set of measures 
$\{ \mu_1\otimes\g\L^p \}_{\g\in E_M}$ is tight by point 4).

As for point 6), we note that the set
$$\{  \g\in\dcal_{\mu_1,\mu_2}\st I(\g)\le M
\}   $$
is weakly compact because of point 5) of this lemma and point 2) of lemma 1.2. Together with point 2) above, this implies the existence of a minimum. The minimizer is unique since $I$ is a strictly convex functional on the convex set 
$\dcal_{\mu_1,\mu_2}$.

\fin

\vskip 1pc
\noindent {\bf Definitions.} $\bullet$) Let $\mu_1,\mu_2\in\Mt$ and let $h>0$; the first equality below is the definition of ${\cal E}^h$, in the second one we recall the definition of $I$ from lemma 1.3. 
$$\Ec^h(\mu_1,\mu_2)=
\min_{\g\in\dcal_{\mu_1,\mu_2}}I(\g)-\log\left(
\frac{1}{2\pi h}
\right)^\frac{p}{2}  =  $$
$$\min_{\g\in\dcal_{\mu_1,\mu_2}}
\int_{\T^p\times\R^p}\left[
\cinh{v}+\log\g(x,v)
\right]   \g(x,v)\dr\mu_1(x)\dr v-\log\left(
\frac{1}{2\pi h}
\right)^\frac{p}{2}  .  \eqno (1.7)$$
Conventionally, we shall say that the minimum is $+\infty$ if 
$\dcal_{\mu_1,\mu_2}$ is empty; if it is not, the minimum in (1.7) is justified by point 6) of lemma 1.3. Let $\g$ be in the class 
$Trace(a,\d)$ defined in the appendix; formula (A.1) of the appendix implies the first inequality below, while the second one follows from the definition of $T(a,\d)$ and the third one from point 1) of lemma A.3.
$$\int_{\R^p}
A_h(\g,v)\dr v-\log\left(
\frac{1}{2\pi h}
\right)^\frac{p}{2}  \ge
T(a,\d)-\log\left(
\frac{1}{2\pi h}
\right)^\frac{p}{2}\ge
T(a,\d)-T(a,1)\ge 0  .  $$
Together with Fubini and (1.7), this implies that
$$\Ec^h(\mu_1,\mu_2)\ge 0\qquad
\forall\mu_1,\mu_2\in\Mt   .   \eqno (1.8)$$
In the inequality above, an explicit calculation shows that $0$ is reached when 
$\mu_2=\mu_1\ast N(0,h Id)$ where  
$$N(0,h Id)(v)=\left(
\frac{1}{2\pi h}
\right)^\frac{p}{2}    
e^{-\cinh{v}}   .   $$

\noindent $\bullet$ As noted by one of the referees, one can express the cost using relative entropy with respect to the normal distribution $N(0,h Id)$; indeed, it is easy to see that, if 
$\g(x,v)=N(0,h Id)(v)\cdot \r(x,v)$, then
$$I(\g)-\left(
\frac{1}{2\pi h}
\right)^\frac{p}{2}  =
\int_{\T^p\times\R^p}
[\r(x,v)\log\r(x,v)] N(0,h Id)(v)\dr\mu_1(x)\dr v  .  $$
See also the connection with the Feynman-Kac formula in section 5 of [6]. 

\noindent $\bullet$) As in point 6) of lemma 1.3, we shall call 
$\g_{\mu_1,\mu_2}$ the unique $\g$ on which 
$\Ec^h(\mu_1,\mu_2)$ is attained.

\noindent $\bullet$) Let us suppose that 
$\Ec^h(\mu_1,\mu_2)<+\infty$; then 
$1\in L^1(\mu_1\otimes\g_{\mu_1,\mu_2}(x,\cdot)\L^p)$ by the definition of $\dcal_{\mu_1,\mu_2}$, while 
$|v|^2\in L^1(\mu_1\otimes\g_{\mu_1,\mu_2}(x,\cdot)\L^p)$ by (1.4). By H\"older's inequality, this implies that 
$v\in L^1(\mu\otimes\g_{\mu_1,\mu_2}(x,\cdot)\L^p)$; following [15],  we define the $h$-forward velocity $v^h_{\mu_1,\mu_2}$ as
$$v^h_{\mu_1,\mu_2}(x)=\frac{1}{h}
\int_{\R^p}v\g_{\mu_1,\mu_2}(x,v)\dr v   .   \eqno (1.9)$$

\noindent $\bullet$) Since 
$|v|^2\in L^1(\mu\otimes\g_{\mu_1,\mu_2}(x,\cdot)\L^p)$, we can define the $h$-covariance matrix $D^{h,\mu_1,\mu_2}$ whose entries 
$d_{i,j}^{h,\mu_1,\mu_2}$ are given by 
$$d_{i,j}^{h,\mu_1,\mu_2}(x)=\frac{1}{h}
\int_{\R^p}[v_i-(v^h_{\mu_1,\mu_2}(x))_ih]\cdot
[v_j-(v^h_{\mu_1,\mu_2}(x))_jh]\g_{\mu_1,\mu_2}(x,v)\dr v  .  $$
If $\mu_t$ is a Borel curve in $\Mt$, we define its $h$-forward velocity as 
$$v^h(t,x)=v^h_{\mu_t,\mu_{t+h}}(x)$$
and its $h$-covariance matrix as
$$D^h(t,x)=D^{h,\mu_t,\mu_{t+h}}(x)   .  $$

\vskip 2pc
\centerline{\bf \S 2}
\centerline{\bf Measurability}
\vskip 1pc

We want to prove the following proposition.

\prop{2.1} Let $\fun{\mu}{(a,b)}{\Mt}$ be a Borel curve of measures and let us suppose that 
$\Ec^h(\mu_t,\mu_{t+h})<+\infty$ for $\L^1$ a. e. $t\in(a,b-h)$. Then, the maps $\fun{}{(t,x)}{v^h(t,x)}$ and 
$\fun{}{(t,x)}{D^h(t,x)}$ (which we defined at the end of the last section) are Borel, up to redefining them on a set of null 
$\L^1\otimes\mu_t$ measure. 

\rm

We shall need the following lemma. 

\lem{2.2} Let $h>0$ and let $v^h_{\mu_1,\mu_2}$, 
$D^{h,\mu_1,\mu_2}$ be defined as at the end of the last section. Let us set 
$$E=\{
(\mu_1,\mu_2)\in\Mt\times\Mt
\st \Ec^h(\mu_1,\mu_2)<+\infty
\}   .  $$ 
Then, the three maps
$$\fun{\Phi_{vel}}{E}{\R}  ,\qquad
\fun{\Phi_{vel}}{(\mu_1,\mu_2)}{
\int_{\T^p}v^h_{\mu_1,\mu_2}(x)\dr\mu_1(x)
} , $$
$$\fun{\Phi^2_{vel}}{E}{\R}  ,\qquad
\fun{\Phi^2_{vel}}{(\mu_1,\mu_2)}{
\int_{\T^p}|v^h_{\mu_1,\mu_2}(x)|^2\dr\mu_1(x)
} , $$
$$\fun{\Phi_{cov}}{E}{\R^{n^2}}  ,\qquad
\fun{\Phi_{cov}}{(\mu_1,\mu_2)}{
\int_{\T^p}D^{h,\mu_1,\mu_2}(x)\dr\mu_1(x)
}  $$
are Borel, while the map 
$$\fun{\Ec^h}{\Mt\times\Mt}{[0,+\infty)}  $$
is lower semicontinuous.

\rm

\vskip 1pc

The idea of the proof is the following. In lemma 2.3 below by 
$\inf$-convolution we are going to find Lipschitz functions 
$c_\l(\mu_1,\mu_2)$ such that, for all couples $(\mu_1,\mu_2)$,
$c_\l(\mu_1,\mu_2)\nearrow\Ec^h(\mu_1,\mu_2)$ as 
$\l\nearrow+\infty$. Clearly, this will yield that 
$\Ec^h$ is l. s. c.. Next, we are going to see that the function 
$\g^\l_{\mu_1,\mu_2}$ on which $c_\l(\mu_1,\mu_2)$ is attained depends continuously on 
$(\mu_1,\mu_2)$ (lemma 2.5) and that it converges to the minimizer of $\Ec^h(\mu_1,\mu_2)$ as $\l\nearrow+\infty$; this will imply that the minima depend in a Borel way on the parameters 
$(\mu_1,\mu_2)$ (corollary 2.6). We begin with a few definitions.

\vskip 1pc

\noindent{\bf Definitions.} $\bullet$) This was introduced in [3], where it is called "push-forward by plans". Let 
$\mu_1,\tilde\mu_1\in\Mt$ and let $\G$ be a transfer plan from  
$\mu_1$ to $\tilde\mu_1$, minimal for the 1-Wasserstein distance $d_1$; we disintegrate 
$\G$ as $\G=\G_y\otimes\tilde\mu_1$. In the following, we shall  reserve the variables $x$ and $y$ for integration in $\mu_1$ and 
$\tilde\mu_1$ respectively. Let $\g\in\dcal_{\mu_1}$; we define
$$\tilde\g(y,v)=\int_{\T^p}\g(x,v)\dr\G_y(x)  .  \eqno (2.1)$$
This is just a generalized way of composing with a map: indeed, if $\G$ is induced by an invertible map $g$, then 
$\tilde\g(y,v)=\g(g^{-1}(y),v)$. 

We refer the reader to [6] for the easy proof that 
$\tilde\g\in\dcal_{\tilde\mu_1}$.

\noindent $\bullet$) We define $Den^\prime$ as the set of all the Borel functions $\fun{\g}{\T^p\times\R^p}{[0,+\infty)}$ such that 
$\g(x,\cdot)$ is a probability density on $\R^p$ for all $x\in\T^p$.

\noindent $\bullet$) For $\l>0$ we define the map
$$\fun{\Psi_\l}{\Mt\times\Mt\times Den^\prime}{\R}$$
$$\Psi_\l(\mu_1,\mu_2,\g)=
\int_{\T^p\times\R^p}A_h(\g,(x,v))\dr\mu_1(x)\dr v+
\l d_1(\mu_1\ast\g,\mu_2)-
\left(
\frac{1}{2\pi h}
\right)^\frac{p}{2}   .   $$

\noindent $\bullet$) We define
$$c_\l(\mu_1,\mu_2)=\inf_{\g\in Den^\prime} 
\Psi_\l(\mu_1,\mu_2,\g)  .  $$

\lem{2.3} 1) The $\inf$ in the definition of $c_\l$ is attained on a function $\g^\l_{\mu_1,\mu_2}\in Den^\prime$; this function is unique up to 
$\mu_1\otimes\L^p$-null sets.

\noindent 2) Let $\g\in \dcal_{\mu_1}$ and let $\tilde\g$ be defined as in (2.1); then, $\tilde\g\in \dcal_{\tilde\mu_1}$ and there is a constant $L=L(\l)>0$ such that
$$\Psi_\l(\tilde\mu_1,\mu_2,\tilde\g)\le\Psi_\l(\mu_1,\mu_2,\g)+
L(\l)d_1(\mu_1,\tilde\mu_1)  .  $$

\noindent 3) The function $c_\l$ is Lipschitz in both arguments for the 1-Wasserstein distance; the Lipschitz constant is 
$\max(L(\l),\l)$.
$$c_\l(\mu_1,\mu_2)\le\Ec^h(\mu_1,\mu_2)
\qquad\forall\mu_1,\mu_2\in\Mt   .  \leqno 4)$$
$$\lim_{\l\tends+\infty}c_\l(\mu_1,\mu_2)=\Ec^h(\mu_1,\mu_2)
\qquad\forall\mu_1,\mu_2\in\Mt . 
\leqno 5)  $$

\noindent 6) Let $\Ec^h(\mu_1,\mu_2)<+\infty$ and let 
$\g_{\mu_1,\mu_2}$ be as in point 6) of lemma 1.3; then, 
$$\mu_1\otimes\g^\l_{\mu_1,\mu_2}\L^p\tends
\mu_1\otimes\g_{\mu_1,\mu_2}\L^p  \txt{narrowly as} 
\l\tends+\infty . 
$$

\noindent 7) For any fixed $\l>0$, 
$$\sup\{
||\g^\l_{\mu_1,\mu_2}||_{L^\infty(\mu_1\otimes\L^p)}
\st \mu_1,\mu_2\in\Mt
\}   <+\infty   .  $$

\proof We note that $\fun{}{\g}{\Psi_\l(\mu_1,\mu_2,\g)}$ has the form
$$\fun{}{\g}{
\int_{\T^p\times\R^p} A_h(\g,(x,v))\dr\mu_1(x)\dr v+U(\mu_1\ast\g)
}    \eqno (2.2)$$
and that $U$ is a Lipschitz function for the 1-Wasserstein distance; namely, for $\mu_2$ fixed, 
$$U(\mu)=\l d_1(\mu,\mu_2)-
\left(
\frac{1}{2\pi h}
\right)^\frac{p}{2}   .   $$
Since $U$ is Lipschitz, proposition 1.4 of [6] holds, yielding  existence in point 1). As for the uniqueness, it suffices to note that the functional of (2.2) is the sum of two terms, the integral and $U$; both are convex in $\g$, the first one strictly. 

Point 2) is proven in proposition 2.3 of [6]. 

As for point 3), from point 2) of this lemma it is easy to deduce (see [6] for the complete argument) that $c_\l$ is 
$L(\l)$-Lipschitz in the first variable. It is $\l$-Lipschitz in the second one because of the special form of the final condition $U$.

We prove point 4). Let $\g_{\mu_1,\mu_2}$ minimize in the definition of $\Ec^h(\mu_1,\mu_2)$; the inequality below comes from the definition of $c_\l(\mu_1,\mu_2)$ as an $\inf$, the first equality from the fact that $\g_{\mu_1,\mu_2}\in\dcal_{\mu_1,\mu_2}$, the second one from the fact that $\Ec^h(\mu_1,\mu_2)$ is attained on 
$\g_{\mu_1,\mu_2}$.
$$c_\l(\mu_1,\mu_2)\le
\int_{\T^p\times\R^p} A_h(\g_{\mu_1,\mu_2},(x,v))
\dr\mu_1(x)\dr v+
\l d_1(\mu_1\ast\g_{\mu_1,\mu_2},\mu_2)-\left(
\frac{1}{2\pi h}
\right)^\frac{p}{2}=$$
$$\int_{\T^p\times\R^p} A_h(\g_{\mu_1,\mu_2},(x,v))
\dr\mu_1(x)\dr v-\left(
\frac{1}{2\pi h}
\right)^\frac{p}{2}= \Ec^h(\mu_1,\mu_2)   . $$
Having thus proven point 4), point 5) reduces to show that
$$\liminf_{\l\tends+\infty}c_\l(\mu_1,\mu_2)\ge
\Ec^h(\mu_1,\mu_2) . $$
Let us suppose by contradiction that this is not the case; in other words, there are $\e>0$ (or $M>0$), a sequence $\l_n\nearrow+\infty$ and minima $\g^{\l_n}_{\mu_1,\mu_2}$ of 
$\Psi_{\l_n}(\mu_1,\mu_2,\cdot)$ such that
$$\int_{\T^p\times\R^p}A_h(\g^{\l_n}_{\mu_1,\mu_2},(x,v))
\dr\mu_1(x)\dr v+
\l_n d_1(\mu_1\ast\g^{\l_n}_{\mu_1,\mu_2},\mu_2)-
\left(\frac{1}{2\pi h}\right)^\frac{p}{2}\le$$
$$\left\{
\eqalign{
\Ec^h(\mu_1,\mu_2)-\e &\txt{if}\Ec^h(\mu_1,\mu_2)<+\infty\cr
M&\txt{otherwise.}
}
\right.      \eqno (2.3)  $$
Now we can apply point 5) of lemma 1.3 and get that, up to subsequences, 
$\g^{\l_n}_{\mu_1,\mu_2}\weak\g$ in $L^1(\mu_1\otimes\L^p)$; by point 2) of lemma 1.3, 
$$\int_{\T^p\times\R^p}A_h(\g,(x,v))\dr\mu_1(x)\dr v-
\left(
\frac{1}{2\pi h}
\right)^\frac{p}{2}\le
\left\{
\eqalign{
\Ec^h(\mu_1,\mu_2)-\e &\txt{if}\Ec^h(\mu_1,\mu_2)<+\infty\cr
M&\txt{otherwise.}
}
\right.       \eqno (2.4)  $$
Since $\l_n\nearrow+\infty$, (1.5) and (2.3) imply that 
$d_1(\mu_1\ast\g_n,\mu_2)\tends 0$; by point 1) of lemma 1.2, we get that $\mu_1\ast\g=\mu_2$. Thus, $\g\in\dcal_{\mu_1,\mu_2}$ and satisfies (2.4): we have reached a contradiction with the definition of $\Ec^h(\mu_1,\mu_2)$ as a minimum. 

We prove point 6). Since narrow convergence is metric (see for instance [2], remark 5.1.1), it suffices to prove that, for any 
$\l_n\nearrow+\infty$ there is a subsequence $\l_{n^\prime}$ such that 
$\mu_1\otimes\g^{\l_{n^\prime}}_{\mu_1,\mu_2}\L^p\tends\mu_1\otimes\g_{\mu_1,\mu_2}\L^p$ narrowly.  Let $\l_n\tends+\infty$; using the fact that 
$\Ec^h(\mu_1,\mu_2)<+\infty$, we can see as in the proof of point 5) that, for a subsequence $\{ n^\prime \}$, 
$\g^{\l_{n^\prime}}_{\mu_1,\mu_2}\weak\g$ in 
$L^1(\mu_1\otimes\L^p)$, that $\g\in\dcal_{\mu_1,\mu_2}$ and that
$$\int_{\T^p\times\R^p}A_h(\g,(x,v))\dr\mu_1(x)\dr v-
\left(
\frac{1}{2\pi h}
\right)^\frac{p}{2}\le\Ec^h(\mu_1,\mu_2)  .  $$
By the uniqueness of point 6) of lemma 1.3, we get that 
$\g=\g_{\mu_1,\mu_2}$. Thus, 
$\g^{\l_{n^\prime}}_{\mu_1,\mu_2}\weak\g_{\mu_1,\mu_2}$ in 
$L^1(\mu_1\otimes\L^p)$; in particular, if 
$\fun{f}{\T^p\times\R^p}{\R}$ is a bounded continuos function, we see that
$$\int_{\T^p\times\R^p}
f(x,v)\g^{\l_{n^\prime}}_{\mu_1,\mu_2}(x,v)\dr\mu_1(x)\dr v\tends
\int_{\T^p\times\R^p}
f(x,v)\g_{\mu_1,\mu_2}(x,v)\dr\mu_1(x)\dr v$$
implying point 6).

As for point 7), this is proposition 2.8 of [6]. 

\fin

Points 3), 4) and 5) of the last lemma imply the last assertion of proposition 2.1; we state it as a separate corollary. 

\cor{2.4} The function
$$\fun{\Ec^h}{\Mt\times\Mt}{[0,+\infty)}$$
is lower semicontinuous.

\rm

\lem{2.5} Let $\l>0$ be fixed and let $\g^\l_{\mu_1,\mu_2}$ be the unique minimizer in the definition of $c_\l(\mu_1,\mu_2)$. Let 
${\cal M}_1(\T^p\times\R^p)$ be the space of the Borel probability measures on $\T^p\times\R^p$ with the topology of  narrow convergence. Then, the map from $\Mt\times\Mt$ to 
${\cal M}_1(\T^p\times\R^p)$ given by 
$$\fun{b^\l}{(\mu_1,\mu_2)}{
\mu_1\otimes\g^\l_{\mu_1,\mu_2}(x,\cdot)\L^p
}    $$
is continuous.

\proof Let $\mu_1^n\tends\mu_1$, $\mu_2^n\tends\mu_2$; we must prove that the sequence 
$\mu_1^n\otimes\g^\l_{\mu_1^n,\mu_2^n}(x,\cdot)\L^p$ converges narrowly to 
$\mu_1\otimes\g^\l_{\mu_1,\mu_2}(x,\cdot)\L^p$. Since $c^\l$ is continuous, we can as well suppose that 
$c^\l(\mu_1^n,\mu^n_2)\le M^\prime$ for some $M^\prime>0$. Together with the definition of $\Psi_\l$ this implies that 
$$\int_{\T^p\times\R^p}
A_h(\g^\l_{\mu_1^n,\mu_2^n},(x,v))
\dr\mu_1^n(x)\dr v\le M\qquad\forall n\ge 1  .  \eqno (2.5)$$
Let $\G^n$ be a plan from $\mu_1^n$ to $\mu_1$ optimal for the 1-Wasserstein distance; let us disintegrate it as 
$\G^n=\G^n_y\otimes\mu_1$ and let us define 
$\tilde\g_n\in{\dcal_{\mu_1}}$ as in (2.1), i. e. 
$$\tilde\g_n(y,v)=
\int_{\T^p}\g^\l_{\mu_1^n,\mu_2^n}(x,v)\dr\G^n_y(x)  .  $$

\noindent{\bf Step 1.} We begin to prove that 
$\mu_1\otimes\tilde\g_n(x,\cdot)\L^p$ converges narrowly to 
$\mu_1\otimes\g^\l_{\mu_1,\mu_2}(x,\cdot)\L^p$.

The inequality below follows by point 2) of lemma 2.3 and by the special form of the final condition. Point 3) of lemma 2.3 implies the limit, while the equality follows because 
$\g^\l_{\mu_1^n,\mu_2^n}$ is minimal. 
$$\Psi_\l(\mu_1,\mu_2,\tilde\g_n)\le
\Psi_\l(\mu_1^n,\mu_2^n,\g^\l_{\mu_1^n,\mu_2^n})+
L(\l)d_1(\mu_1^n,\mu_1)+\l d_1(\mu_2^n,\mu_2) = $$
$$c_\l(\mu_1^n,\mu_2^n)+
L(\l)d_1(\mu_1^n,\mu_1)+\l d_1(\mu_2^n,\mu_2)\tends
c_\l(\mu_1,\mu_2)  .  $$
In other words, $\tilde\g_n$ is a minimizing sequence for 
$c_\l(\mu_1,\mu_2)$; this implies (proposition 1.4 of [6]) that 
$\tilde\g_n\weak\g^\l_{\mu_1,\mu_2}$ in $L^1(\mu_1\otimes\L^p)$. As we have seen at the end of the proof of lemma 2.3, this implies that $\mu_1\otimes\tilde\g_n(x,\cdot)\L^p$ converges narrowly to 
$\mu_1\otimes\g^\l_{\mu_1,\mu_2}(x,\cdot)\L^p$.

\noindent{\bf Step 2.} By step 1 it suffices to prove that, if 
$\fun{f}{\T^p\times\R^p}{\R}$ is a bounded continuous function, then 
$$\int_{\T^p\times\R^p}
f(x,v)\g^\l_{\mu_1^n,\mu_2^n}(x,v)\dr\mu_1^n(x)\dr v-
\int_{\T^p\times\R^p}
f(y,v)\tilde\g_n(x,v)\dr\mu_1(y)\dr v\tends 0  .  
\eqno (2.6)  $$
Note that $\mu_1^n\otimes\g^\l_{\mu_1^n,\mu_2^n}\L^p$ is tight by (1.4) and (2.5); using this, we easily see that it suffices to consider a uniformly continuous $f$. 

The first equality below is the definition of $\tilde\g_n$, the second one comes from the fact that $\G^n=\G^n_y\otimes\mu_1$ and the fact that the first marginal of $\G^n$ is $\mu^n_1$. 
$$\int_{\T^p\times\R^p}
f(x,v)\g^\l_{\mu_1^n,\mu_2^n}(x,v)\dr\mu_1^n(x)\dr v-
\int_{\T^p\times\R^p}
f(y,v)\tilde\g_n(y,v)\dr\mu_1(y)\dr v=$$
$$\int_{\T^p\times\R^p}
f(x,v)\g^\l_{\mu_1^n,\mu_2^n}(x,v)\dr\mu_1^n(x)\dr v-
\int_{\T^p\times\R^p}
f(y,v)\dr\mu_1(y)\dr v
\int_{\T^p}\g^\l_{\mu_1^n,\mu_2^n}(x,v)\dr\G^n_y(x)=$$
$$\int_{\T^p\times\T^p\times\R^p}[
f(x,v)-f(y,v)
]   \g^\l_{\mu_1^n,\mu_2^n}(x,v)
\dr\G^n(x,y)\dr v      .   \eqno (2.7)$$
Since $f$ is uniformly continuous, for all $\e>0$ we can find $\d>0$ such that
$$\txt{if}|x-y|<\d
\txt{then} ||f(x,\cdot)-f(y,\cdot)||_{C^0(\R^p)}\le\e  .  \eqno (2.8)$$
We set
$$A_\d=\{
(x,y)\in\T^p\times\T^p\st |x-y|_{\T^p}<\d
\}   .     $$
Formula (2.7) implies the first inequality below; the second one follows by (2.8), the fact that $f$ is bounded and the fact that 
$\g^\l_{\mu_1^n,\mu_2^n}(x,\cdot)$ is a probability density. 
$$\left\vert
\int_{\T^p\times\R^p}
f(x,v)\g^\l_{\mu_1^n,\mu_2^n}(x,v)\dr\mu_1^n(x)\dr v-
\int_{\T^p\times\R^p}
f(y,v)\tilde\g_n(y,v)\dr\mu_1(y)\dr v
\right\vert   \le$$
$$\int_{A_\d\times\R^p}
|f(x,v)-f(y,v)|\g^\l_{\mu_1^n,\mu_2^n}(x,v)
\dr\G^n(x,y)\dr v+
\int_{A_\d^c\times\R^p}
|f(x,v)-f(y,v)|\g^\l_{\mu_1^n,\mu_2^n}(x,v)
\dr\G^n(x,y)\dr v   \le$$
$$\int_{A_\d}\e\dr\G^n(x,y)+
2||f||_\infty\int_{A_\d^c}\dr\G^n(x,y)  .  $$
In other words, (2.6) follows if we prove that 
$\G^n(A_\d^c)\tends 0$; but this comes from the Chebyshev inequality below. 
$$\d\G^n(A_\d^c)\le
\int_{A_\d^c}|x-y|_{\T^p}\dr\G^n(x,y)\le$$
$$\int_{\T^p\times\T^p}
|x-y|_{\T^p}\dr\G^n(x,y)=
d_1(\mu_1^n,\mu_1)\tends 0  .   $$

\fin

\cor{2.6} Let 
$$E=\{
(\mu_1,\mu_2)\in\Mt\times\Mt\st\Ec^h(\mu_1,\mu_2)<+\infty
\}  .  $$
Then, the map
$$\fun{b}{E}{
{\cal M}_1(\T^p\times\R^p)
}  \qquad
\fun{b}{(\mu_1,\mu_2)}{\mu_1
\otimes\g_{\mu_1,\mu_2}(x,\cdot)\L^p}  $$
is Borel; we have endowed $\Mt$ with the weak$\ast$ and 
${\cal M}_1(\T^p\times\R^p)$  with the narrow topology. 

\proof By lemma 2.5, the map $b^\l$ is continuous. Thus, it suffices to show that, for all $\mu_1,\mu_2$ such that 
$\Ec^h(\mu_1,\mu_2)<+\infty$, we have that
$$\mu_1\otimes\g^\l_{\mu_1,\mu_2}(x,\cdot)\L^p\tends
\mu_1\otimes\g_{\mu_1,\mu_2}(x,\cdot)\L^p  $$
narrowly as $\l\tends+\infty$. But this is the content of point 6) of lemma 2.3. 

\fin

\noindent{\bf Proof of lemma 2.2.} The map 
$\Ec^h $ is l. s. c. by corollary 2.4. We want to prove that 
$\Phi_{vel}$, $\Phi_{vel}^2$ and $\Phi_{cov}$ are Borel. We prove that $\Phi^2_{vel}$ is Borel, since the other cases are analogous. This follows easily by corollary 2.6 and the fact that the map
$$\colon
\mu_1\otimes\g_{\mu_1,\mu_2}(x,\cdot)\L^p
\rightarrow
\int_{\T^p\times\R^p}|v|^2\g_{\mu_1,\mu_2}(x,v)\dr v\dr\mu_1(x)$$
is l. s. c..

\fin

\noindent{\bf Proof of proposition 2.1.} Let $\g_t$ be the unique (up to $\mu_t$-null sets) minimizer in the definition of 
$\Ec^h(\mu_t,\mu_{t+h})$; we are going to prove that $\g_t$ has a Borel version. 

By corollary 2.6, we can define a measure $\d$ on 
$(a,b)\times\T^p\times\R^p$ by
$$\int_{(a,b)\times\T^p\times\R^p}f(t,x,v)\dr\d(t,x,v)=
\int_a^b\dr t\int_{\T^p\times\R^p}f(t,x,v)\g_t(x,v)\dr\mu_t(x)\dr v$$
for all continuous, compactly supported functions $f$. Cearly, 
$\d<<\L^1\otimes\mu_t\otimes\L^p$; by the Radon-Nikodym theorem, we can write 
$\d=\L^1\otimes\mu_t\otimes\tilde\g_t(x,v)\L^p$; now $\tilde\g_t$ is a Borel version of $\g_t$. 

For simplicity, from now on we shall drop the tilde from $\tilde\g_t$, as if $\g_t$ were already Borel. We prove that $v^h(t,x)$ is Borel, up to modifying it on a set of null $\L^1\otimes\mu_t$ measure. We shall forego the proof that $D^h(t,x)$ is Borel, which is similar. 

First of all, it suffices to find Borel sets $A_n\subset(a,b)$ such that 

\noindent 1) $\L^1((a,b)\setminus A_n)\tends 0$ and 

\noindent 2) $v^h$ is Borel on $A_n\times\T^p$.

\noindent Since 
$$v^h(t,x)=\frac{1}{h}\int_{\R^p}v\g_t(x,v)\dr v  ,$$
Fubini's theorem implies that (up to redefining it on a set of null 
$\L^1\otimes\mu_t$ measure) $v^h$ is Borel on  $A_n\times\T^p$ if 
$v\in L^1(A_n\times\T^p\times\R^p,\L^1\otimes\mu_t\otimes
\g_t(x,\cdot)\L^p)$; actually, we shall show that 
$v\in L^2(A_n\times\T^p\times\R^p,\L^1\otimes\mu_t\otimes
\g_t(x,\cdot)\L^p)$. 

We set 
$$A_n=\{
t\in(a,b)\st \Ec^h(\mu_t,\mu_{t+h})\le n
\}   .   $$
This is a Borel set by corollary 2.4; the sets $A_n$ invade 
$(a,b)$ since $\Ec^h(\mu_t,\mu_{t+h})<+\infty$ for a. e. 
$t\in(a,b)$ by hypothesis. This yields point 1) above; point 2) follows since point 4) of lemma 1.3 implies the first inequality below.
$$\int_{A_n\times\T^p\times\R^p}
\cinh{v}\g_t(x,v)\dr\mu_t(x)\dr v\le
\int_{A_n}B(n)\dr t\le B(n)(b-a)  .  $$

\fin

\vskip 1pc
\centerline{\bf \S 3}
\centerline{\bf Forward velocity and diffusion matrix}
\vskip 1pc

In this section, we define the cost of a curve of measures $\mu_t$; in proposition 3.1 below we shall see that, if the cost of 
$\mu_t$ is finite, the $h$-forward velocity we defined in section 1 is bounded in $L^2$; taking limits, we shall get an instantaneous forward velocity. A similar argument will yield (proposition 3.2 below) that the $h$-diffusion matrix converges to the identity. 

\vskip 1pc

\noindent{\bf Definition.} Let $\fun{\mu}{(a,b)}{\Mt}$ be Borel and let $h\in(0,b-a)$. We define
$$\Ec^h_{(a,b-h)}(\mu)=\frac{1}{h}\int_a^{b-h}
\Ec^h(\mu_t,\mu_{t+h})\dr t$$
where the function $\Ec^h(\mu_t,\mu_{t+h})$ has been defined in (1.7). We note that this integral is well defined, though possibly 
$+\infty$, because the map $\fun{}{t}{\Ec^h(\mu_t,\mu_{t+h})}$ non negative by (1.8); it is Borel because it is the composition of 
$\Ec^h$, which is l. s. c. by lemma 2.2, with the Borel function $\mu_t$.  We also set
$$\Ec_{(a,b)}(\mu)=\liminf_{h\tends 0}
\Ec^h_{(a,b-h)}(\mu)  .  $$

Instead of proving that $\Ec_{(a,b)}$ is lower semicontinuous, we are going to relax it, so that semicontinuity will be automatic. We define
$${\cal C}_{(a,b)}(\mu)=\inf
\liminf_{n\tends+\infty}\Ec_{(a,b)}(\mu_n)  $$
where the $\inf$ is over all sequences 
$\fun{\mu_n}{(a,b)}{\Mt}$ converging uniformly to $\mu$ for the 2-Wasserstein distance $d_2$. 

\vskip 1pc

Our aim is to study the relaxed cost ${\cal C}_{(a,b)}$; to do this, we need some preliminary knowledge on $\Ec_{(a,b)}$; this will take all of this section and some of the next one. 

\prop{3.1} Let $\fun{\mu}{(a,b)}{\Mt}$ be a Borel curve of measures and let the $h$-forward velocity $v^h(t,x)$ be as in the definition at the end of section 1. Then, the following two points hold.

\noindent 1) If $\Ec^h_{(a,b-h)}(\mu)<+\infty$, then $v^h(t,\cdot)$ is defined for $\L^1$ a. e. $t\in(a,b-h)$ and
$$\int_a^{b-h}\dr t\int_{\T^p}
\cin{v^h(t,x)}\dr\mu_t(x)\le\Ec^h_{(a,b-h)}(\mu)  .  
\eqno (3.1)  $$

\noindent 2) Let $\Ec_{(a,b)}(\mu)<+\infty$; then, there is a subsequence $h_n\searrow 0$ and 
$X\in L^2((a,b)\times\T^p,\L^1\otimes\mu_t)$ such that
$$v^{h_n}\weak X\txt{in}
L^2((a,b)\times\T^p,\L^1\otimes\mu_t)   .  $$

\noindent 3) Let $X$ be as in point 2) above. Then, 
$$\int_a^b\dr t\int_{\T^p}
|X(t,x)|^2\dr\mu_t(x)\le \Ec_{(a,b)}(\mu)  .  $$

\proof We begin with point 1); let us show that the integral defining $v^h(t,\cdot)$ converges. We recall that, by definition, 
$$v^h(t,x)=v^h_{\mu_t,\mu_{t+h}}(x)$$
and that the integral defining $v^h_{\mu_t,\mu_{t+h}}$ converges  if 
$\g_{\mu_t,\mu_{t+h}}$ has finite second moment; by (1.4), this is true if $\Ec^h(\mu_t,\mu_{t+h})<+\infty$. Thus, it suffices to show that 
$\Ec^h(\mu_t,\mu_{t+h})<+\infty$ for $\L^1$ a. e. $t\in(a,b-h)$; since ${\cal E}^h(\mu_t,\mu_{t+h})\ge 0$, this follows from the formula below, where the inequality is our hypothesis and the equality is the definition of 
$\Ec^h_{(a,b-h)}(\mu)$. 
$$\frac{1}{h}\int_a^{b-h}\Ec^h(\mu_t,\mu_{t+h})\dr t=
\Ec^h_{(a,b-h)}(\mu)< +\infty .  $$
To prove (3.1), we set 
$$\g^h_t\colon=\g_{\mu_t,\mu_{t+h}}$$
where $\g_{\mu_t,\mu_{t+h}}$ has been defined in point 6) of lemma 1.3. We recall that $\g^h_t(x,\cdot)\in Mean(v^h(t,x))$ by the definition of $v^h(t,x)$ in section 1 and of $Mean(v^h(t,x))$ in the appendix. 
We define the trace of the variance as in the appendix
$$p\d(t,x)=\frac{1}{h}\int_{\R^p}
|v-hv^h(t,x)|^2\g^h_t(x,v)\dr v  .  $$
The inequality below follows by (A.1) of the appendix and the definition of $T(v^h(t,x),\d(x))$; the equality is (1.7).
$$\frac{1}{h}\int_a^{b-h}\Ec^h(\mu_t,\mu_{t+h})\dr t=
\frac{1}{h}\int_a^{b-h}\dr t\int_{\T^p}\dr\mu_t(x)
\left[
\int_{\R^p}A_h(\g^h_t,(x,v))\dr v-
\log\left( \frac{1}{2\pi h} \right)^\frac{p}{2}
\right]   \ge $$
$$\frac{1}{h}\int_a^{b-h}\dr t\int_{\T^p}
\left[
p\frac{\d(t,x)-1}{2}+\frac{h}{2}|v^h(t,x)|^2+
\log\left( \frac{1}{\d(t,x)} \right)^\frac{p}{2}
\right]  \dr\mu_t(x)  .  $$
Since
$$p\frac{\d-1}{2}+\log\left(\frac{1}{\d}\right)^\frac{p}{2}\ge 0
\qquad\forall\d>0  ,  $$
formula (3.1) follows. 

We prove point 2). Since $\Ec_{(a,b)}(\mu)<+\infty$, we can find 
$h_n\searrow 0$ and $M>0$ such that
$$\Ec^{h_n}_{(a,b-h_n)}(\mu)\le M  \qquad\forall n\ge 1 . 
\eqno (3.2) $$
By (3.1), this implies that
$$\int_a^{b-h_n}\dr t\int_{\T^p}
\cin{v^{h_n}(t,x)}\dr\mu_t(x)\le M\qquad  \forall n\ge 1  .  $$
In other words, $v^{h_n}$ is bounded in 
$L^2((a,b)\times\T^p,\L^1\otimes\mu_t)$; thus, it is weakly compact and point 2) follows. 

Point 3) follows immediately from points 1) and 2) and the lower semicontinuity of the $L^2$ norm under weak convergence. 

\fin

This calls for a definition.

\vskip 1pc
\noindent {\bf Definition.} Let $X$ be as in point 2) of proposition 3.1; we say that $X$ is a forward velocity of the curve $\mu$. 

\vskip 1pc

Note that there is no uniqueness for the forward velocity: different sequences $h_n\searrow 0$ may yield different forward velocities.

\prop{3.2} Let $\fun{\mu}{(a,b)}{\Mt}$ be Borel with 
${\cal E}_{(a,b)}(\mu)<+\infty$; let the $h$-correlation matrix 
$D^h(t,x)=(d_{i,j}^h(t,x))_{i,j}$ be defined as in section 1. Let 
$h_n\searrow 0$ be the sequence of proposition 3.1; then,  
$$D^{h_n}\tends Id\txt{in}
L^1((a,b)\times\T^p,\L^1\otimes\mu_t)  .  \eqno (3.3)  $$

\proof We begin to show that any term on the diagonal, say 
$d_{i,i}^{h_n}$, tends to 1. Since $\{ h_n \}$ is the sequence of proposition 3.1, (3.2) holds and this yields the first inequality below; the equality is (1.7); the second inequality comes from (A.4) and the last one comes from the definition of $B_{diag}$ in the appendix.
$$M\ge\frac{1}{h_n}\int_a^{b-h_n} 
\Ec^{h_n}(\mu_t,\mu_{t+h_n})\dr t=$$
$$\frac{1}{h_n}\int_a^{b-h_n}\dr t\int_{\T^p}\dr\mu_t(x)
\left[
\int_{\R^p}A_{h_n}(\g^{h_n}_t(x,\cdot),v)\dr v-
\log\left( \frac{1}{2\pi h_n} \right)^\frac{p}{2}
\right]   \ge  $$
$$\frac{1}{h_n}\int_a^{b-h_n}\dr t\int_{\T^p}
\left[
B_{diag}(v^{h_n}(t,x),d^{h_n}_{i,i}(t,x))-
\log\left( \frac{1}{2\pi h_n} \right)^\frac{p}{2}
\right]    \dr\mu_t(x)  \ge$$
$$\frac{1}{h_n}\int_a^{b-h_n}\dr t\int_{\T^p}
\left[
B_{diag}(v^{h_n}(t,x),d^{h_n}_{i,i}(t,x))-
B_{diag}(v^{h_n}(t,x),1)
\right]    \dr\mu_t(x)   .    \eqno (3.4)  $$
As in lemma A.3, we shall denote by $D_j$ a constant independent of everything. The last formula and (A.11) imply the first inequality below, while the second one comes from H\"older.
$$h_nM\ge D_1\int_a^{b-h_n}\dr t
\int_{\{ x\in\T^p\st d_{i,i}^{h_n}(t,x)\le 2 \}}
|d_{i,i}^{h_n}(t,x)-1|^2\dr\mu_t(x)+$$
$$D_1\int_a^{b-h_n}\dr t
\int_{\{ x\in\T^p\st d_{i,i}^{h_n}(t,x)> 2 \}}
|d_{i,i}^{h_n}(t,x)-1|\dr\mu_t(x)\ge$$
$$\frac{D_1}{b-a}\left[
\int_a^{b-h_n}\dr t
\int_{\{ x\in\T^p\st d_{i,i}^{h_n}(t,x)\le 2 \}}
|d_{i,i}^{h_n}(t,x)-1|\dr\mu_t(x)
\right]^2  +$$
$$D_1\int_a^{b-h_n}\dr t
\int_{\{ x\in\T^p\st d_{i,i}^{h_n}(t,x)> 2 \}}
|d_{i,i}^{h_n}(t,x)-1|\dr\mu_t(x)  .  $$
This clearly implies that $d^{h_n}_{i,i}\tends 1$ in 
$L^1(\L^1\otimes\mu_t)$. Now we tackle the terms off the diagonal.

With the same argument we used for (3.4) we get that, for 
$i\not=j$, 
$$h_nM\ge\int_a^{b-h_n}\dr t\int_{\T^p}
[
B_{off-diag}(v^{h_n}(t,x),d_{i,j}^{h_n}(t,x))-
B_{off-diag}(v^{h_n}(t,x),0)
]   \dr\mu_t(x)  .  $$
By (A.12), this implies the first inequality below, while the second one is H\"older. 
$$h_nM\ge D_1\int_a^{b-h_n}
\int_{\{ x\in\T^p\st d_{i,j}^{h_n}(t,x)\le 1 \}}
|d_{i,j}^{h_n}(t,x)|^2\dr\mu_t(x)+
D_1\int_a^{b-h_n}
\int_{\{ x\in\T^p\st d_{i,j}^{h_n}(t,x)>1 \}}
|d_{i,j}^{h_n}(t,x)|\dr\mu_t(x)   \ge$$
$$\frac{D_1}{b-a}\left[
\int_a^{b-h_n}
\int_{\{ x\in\T^p\st d_{i,j}^{h_n}(t,x)\le 1 \}}
|d_{i,j}^{h_n}(t,x)|\dr\mu_t(x)
\right]^2+
D_1\int_a^{b-h_n}
\int_{\{ x\in\T^p\st d_{i,j}^{h_n}(t,x)>1 \}}
|d_{i,j}^{h_n}(t,x)|\dr\mu_t(x)  .  $$
The formula above implies that $d^{h_n}_{i,j}\tends 0$ in 
$L^1(\L^1\otimes\mu_t)$, and we are done.

\fin

We shall use the estimates above in the Taylor developments of the next section; we shall also need the third-order estimate below.  First of all, we define a function 
$$\fun{l}{\R^p}{\R},\qquad
l(v)=\left\{
\eqalign{
|v|^3&\txt{if}|v|\le 1\cr
1&\txt{if}|v|\ge 1   .  
}
\right.    \eqno (3.5)  $$

\lem{3.3} Let the curve $\mu_t$ and the sequence 
$h_n\searrow 0$ be as in proposition 3.1; then, 
$$\frac{1}{h_n}\int_a^{b-h_n}\dr t\int_{\T^p}\dr\mu_t(x)
\int_{\R^p}l(v)\g^{h_n}_t(x,v)\dr v\tends 0
\txt{as}n\tends+\infty    \eqno (3.6)$$ 
and
$$\frac{1}{h_n}\int_a^{b-h_n}\dr t\int_{\T^p}\dr\mu_t(x)
\int_{B(0,1)^c}|v|^2\g^{h_n}_t(x,v)\dr v\tends 0
\txt{as}n\tends+\infty  .  \eqno (3.7)$$ 

\proof We begin with (3.6); for $v^{h_n}$ defined as in section 1, we define $\d_n(t,x)$ by
$$p\d_n(t,x)=\frac{1}{h_n}\int_{\R^p}|v-h_nv^{h_n}(t,x)|^2
\g^{h_n}_t(x,v)\dr v  .  $$
We are going to split the innermost integral of (3.6) between 
$B(0,r)$ and $B^c(0,r)$; for the integral on $B(0,r)$ we shall use the fact that on this set $l$ is small, if $r$ is small; for the integral on $B^c(0,r)$ we shall prove that the "tail" of $\g^{h_n}$ tends to zero. 

We note that, for $r\in(0,1)$,
$$\frac{1}{h_n}\int_a^{b-h_n}\dr t\int_{\T^p}\dr\mu_t(x)
\int_{\R^p}l(v)\g^{h_n}_t(x,v)\dr v\le$$
$$\frac{1}{h_n}\int_a^{b-h_n}\dr t\int_{\T^p}\dr\mu_t(x)
\int_{B(0,r)}|v|^3\g^{h_n}_t(x,v)\dr v  +
\eqno (3.8)_a$$
$$\frac{1}{h_n}\int_a^{b-h_n}\dr t\int_{\T^p}\dr\mu_t(x)
\int_{B(0,r)^c}\g^{h_n}_t(x,v)\dr v  .
\eqno (3.8)_b$$
The second equality below comes from the the definition of 
$v^{h_n}$, the third one from the definition of $d^{h_n}_{i,i}$.  
$$(3.8)_a\le \frac{r}{h_n}
\int_a^{b-h_n}\dr t\int_{\T^p}\dr\mu_t(x)
\int_{B(0,r)}|v|^2\g^{h_n}_t(x,v)\dr v\le$$
$$\frac{r}{h_n}
\int_a^{b-h_n}\dr t\int_{\T^p}\dr\mu_t(x)
\int_{\R^p}|v|^2\g^{h_n}_t(x,v)\dr v=$$
$$\frac{r}{h_n}
\int_a^{b-h_n}\dr t\int_{\T^p}\dr\mu_t(x)
\int_{\R^p}[
|v-h_n v^{h_n}(t,x)|^2-h_n^2|v^{h_n}(t,x)|^2+
2\inn{h_n v^{h_n(t,x)}}{v}
]    \g^{h_n}_t(x,v)\dr v=$$
$$h_n r\int_a^{b-h_n}\dr t\int_{\T^p}
|v^{h_n}(t,x)|^2\dr\mu_t(x)+
\frac{r}{h_n}\int_a^{b-h_n}\dr t\int_{\T^p}\dr\mu_t(x)
\sum_{i=1}^p\int_{\R^p}|v_i-h_nv^{h_n}_i|^2
\g^{h_n}_t(x,v)\dr v=$$
$$h_n r\int_a^{b-h_n}\dr t\int_{\T^p}
|v^{h_n}(t,x)|^2\dr\mu_t(x)+
r\int_a^{b-h_n}\dr t\int_{\T^p}
\sum_{i=1}^p d_{i,i}^{h_n}(t,x)\dr\mu_t(x)  .  $$
By (3.1) and proposition 3.2, this implies that there is $n_0(r)>0$ such that 
$$(3.8)_a\le 2pr\txt{for} n\ge n_0(r).    \eqno (3.9)$$
We tackle $(3.8)_b$. Let us set
$$\e_n(t,x)=\frac{1}{h_n}\int_{B(0,r)^c}
\g^{h_n}_t(x,v)\dr v  ,\qquad
\d_n(t,x)=\frac{1}{2h_n}\int_{B(0,r)^c}|v|^2
\g^{h_n}_t(x,v)\dr v  .  $$
Clearly,
$$\e_n(t,x)\le\frac{2}{r^2}\d_n(t,x)
\txt{for}
(t,x)\in(a,b)\times\T^p  .  \eqno (3.10)$$
For $\d_0(r,h)$ as in lemma A.4 we define 
$$A^n_t=\{
x\in\T^p\st\e_n(t,x)\in(0,\frac{2}{r^2}\d_0(r,h_n))
\}  ,  $$
$$C^n_t=\{
x\in\T^p\st\e_n(t,x)\ge \frac{2}{r^2}\d_0(r,h_n) 
\}    . $$
By the definition of $\e_n$ this implies that
$$(3.8)_b=\int_a^{b-h_n}\dr t
\int_{
A^n_t\cup C^n_t
}  \e_n(t,x)\dr\mu_t(x)  .  \eqno (3.11)$$
The definition of $A^n_t$ implies the first inequality below; for the limit, we know by lemma A.4 that $\d_0(r,h_n)\tends 0$.
$$\int_a^{b-h_n}\dr t
\int_{A^n_t}\e_n(t,x)\dr\mu_t(x)\le
\frac{2}{r^2}
\int_a^{b-h_n} \d_0(r,h_n)\dr t   \tends 0   .   \eqno (3.12)$$
Note that, if $x\in C^n_t$, then by (3.10) $\d_n(t,x)\ge\d_0(r,h_n)$; by lemma A.4 this implies that
$$\Ec^{h_n}(\mu_t,\mu_{t+h_n})\ge\frac{\d_n(t,x)}{2}  .  $$
This implies the second inequality below, while the first one follows from (3.10).The last inequality below follows from the definition of 
$\Ec^{h_n}_{(a,b-h_n)}(\mu)$, while the limit follows from (3.2). 
$$\int_a^{b-h_n}\dr t\int_{C^n_t}\e_n(t,x)\dr\mu_t(x)\le
\frac{2}{r^2}\int_a^{b-h_n}\dr t\int_{C^n_t}
\d_n(t,x)\dr\mu_t(x)\le$$
$$\frac{4}{r^2}\int_a^{b-h_n}\Ec^{h_n}(\mu_t,\mu_{t+h_n})\dr t\le
\frac{4}{r^2}h_n\Ec^{h_n}_{(a,b-h_n)}(\mu)  \tends 0  .  $$
By (3.11), (3.12) and the last formula we get that 
$$(3.8)_b\tends 0 \txt{as}n\tends+\infty  .  $$
Let $\e>0$ be given; by (3.9) we can find $r>0$ so small that 
$(3.8)_a\le\frac{\e}{2}$; by the last formula, we can choose $n$ so large that $(3.8)_b\le\frac{\e}{2}$; by (3.8), this implies (3.6).

We prove (3.7). We set 
$$\e_n(t,x)=\frac{1}{2h_n}\int_{B(0,1)^c}|v|^2
\g^{h_n}_t(x,v)\dr v   . $$
For $\d_0(1,h)$ as in lemma A.4, we define
$$A^n_t=\{
x\in\T^p\st\e_n(t,x)\le\d_0(1,h_n)
\}  ,  \qquad
C^n_t=\{
x\in\T^p\st\e_n(t,x)>\d_0(1,h_n)
\}  .  $$
These definitions yield the first two equalities below; the inequality follows from lemma A.4. 
$$\frac{1}{h_n} \int_a^{b-h_n}\dr t\int_{\T^p}\dr\mu_t(x)
\int_{B(0,1)^c}|v|^2\g^{h_n}_t(x,v)\dr v=
2\int_a^{b-h_n}\dr t\int_{\T^p}\e_n(t,x)\dr\mu_t(x)=$$
$$2\int_a^{b-h_n}\dr t\int_{A^n_t}\e_n(t,x)\dr\mu_t(x)+
2\int_a^{b-h_n}\dr t\int_{C^n_t}\e_n(t,x)\dr\mu_t(x)\le$$
$$2\d_0(1,h_n)\int_a^{b-h_n}\mu_t(A^n_t)\dr t+
4\int_a^{b-h_n}\dr t
\int_{\T^p}\Ec^{h_n}(\mu_t,\mu_{t+h_n})\dr\mu_t(x).  $$
Now (3.2) implies the inequality below
$$4\int_a^{b-h_n}
\Ec^{h_n}(\mu_t,\mu_{t+h_n})\dr t\le h_nM\tends 0  .  $$
Since $\d_0(1,h_n)\tends 0$ by lemma A.4, the last two formulas imply (3.7). 

\fin

\vskip 2pc
\centerline{\bf \S 4}
\centerline{\bf Curves of finite energy satisfy the Fokker-Planck equation}
\vskip 1pc

We are going to use the results of section 3 to prove the following proposition, which is the direct part of theorem 1.

\prop{4.1} Let $\fun{\mu}{(a,b)}{\Mt}$ be Borel and let 
${\cal C}_{(a,b)}(\mu)<+\infty$; then there is a vector field 
$X\in L^2((a,b)\times\T^p,\L^1\otimes\mu_t)$ such that $\mu_t$ is a weak solution of the Fokker-Planck equation with drift $X$. In other words,
$$\int_a^b\dr t\int_{\T^p}
\partial_t\phi(t,x)\dr\mu_t(x)=$$
$$-\2\int_a^b\dr t\int_{\T^p} \D\phi(t,x)\dr\mu_t(x)-
\int_a^b\dr t\int_{\T^p}\inn{\nabla\phi(t,x)}{X(t,x)}\dr\mu_t(x)
\quad\forall\phi\in C^\infty_c((a,b)\times\T^p) .  \eqno (4.1)$$
Moreover, we have that
$$\int_a^b\dr t\int_{\T^p}|X(t,x)|^2\dr\mu_t(x)\le
{\cal C}_{(a,b)}(\mu)    .    \eqno (4.2)$$

\proof {\bf Step 1.} We begin with the (apparently) stronger hypothesis $\Ec_{(a,b)}<+\infty$; we shall come to the case 
${\cal C}_{(a,b)}<+\infty$ in step 2 below. 

Let $\phi\in C^\infty_0((a,b)\times\T^p)$, let the sequence 
$h_n\searrow 0$ be as in proposition 3.2 and let us set 
$\g^{h_n}_t=\g^{h_n}_{\mu_t,\mu_{t+h_n}}$. The first equality below comes from dominated convergence, the third one from the fact that 
$\g^{h_n}_t\in\dcal_{\mu_t,\mu_{t+h_n}}$ and the last one from the definition of $\mu_t\ast\g^{h_n}_t$.
$$\int_a^b\dr t\int_{\T^p}\partial_t\phi(t,x)\dr\mu_t(x)=
-\lim_{n\tends+\infty}\int_a^b\dr t\int_{\T^p}
\frac{\phi(t-h_n,x)-\phi(t,x)}{h_n}\dr\mu_t(x)=$$
$$-\lim_{n\tends+\infty}\int_a^b\dr t\int_{\T^p}
\phi(t,x)\dr\left(
\frac{\mu_{t+h_n}(x)-\mu_t(x)}{h_n}
\right)   =   $$
$$-\lim_{n\tends+\infty}\int_{a}^b\dr t\int_{\T^p}
\phi(t,x)\dr\left(
\frac{(\mu_{t}\ast\g^{h_n}_t)(x)-\mu_t(x)}{h_n}
\right) =$$
$$-\lim_{n\tends+\infty}\int_a^b\dr t\int_{\T^p}\dr\mu_t(x)
\int_{\R^p}\frac{\phi(t,x+v)-\phi(t,x)}{h_n}
\g^{h_n}_t(x,v)\dr v   .   $$
Thus, (4.1) follows if we prove that 
$$\lim_{n\tends+\infty}\int_a^b\dr t\int_{\T^p}\dr\mu_t(x)
\int_{\R^p}\frac{\phi(t,x+v)-\phi(t,x)}{h_n}
\g^{h_n}_t(x,v)\dr v=$$
$$\int_a^b\dr t\int_{\T^p}[
\2\D\phi+\inn{X}{\nabla\phi}
]  \dr\mu_t(x)   .    \eqno (4.3)$$
We show (4.3). Let the function $l$ be as in (3.5); by a Taylor development, we get that there is $D_2>0$ for which the  inequality below holds; the limit at the end follows from (3.6) and (3.7) of lemma 3.3.
$$\Bigg\vert
\int_a^b\dr t\int_{\T^p}\dr\mu_t(x)\int_{B(0,1)\cup B(0,1)^c}
\Bigg[
\frac{\phi(t,x+v)-\phi(t,x)}{h_n}-  $$
$$\frac{1}{h_n}\inn{\nabla\phi(t,x)}{v}-
\frac{1}{2h_n}\phi^\pprime(t,x)(v,v)
\Bigg]   \g^{h_n}_t(x,v)\dr v
\Bigg\vert   \le  $$
$$D_2\int_a^b\dr t\int_{\T^p}\dr\mu_t(x)\int_{B(0,1)}
\frac{1}{h_n}l(v)\g^{h_n}_t(x,v)\dr v+
D_2\int_a^b\dr t\int_{\T^p}\dr\mu_t(x)\int_{B(0,1)^c}
\frac{1}{h_n}|v|^2\g^{h_n}_t(x,v)\dr v\tends 0  .  $$
Thus, (4.3) follows if we show that
$$\lim_{n\tends+\infty}
\frac{1}{h_n}\int_a^b\dr t\int_{\T^p}\dr\mu_t(x)
\int_{\R^p}[
\inn{\nabla\phi(t,x)}{v}+\2\phi^\pprime(t,x)(v,v)
]  \g^{h_n}_t(x,v)\dr v =$$
$$\int_a^b\dr t\int_{\T^p}[
\2\D\phi+\inn{X}{\nabla\phi}
]   \dr\mu_t(x)   .   \eqno (4.4)$$
We begin with the gradient term; the equality below comes from the definition of $v^{h_n}$ in section 1, while the limit comes from point 2) of proposition 3.1.
$$\int_a^b\dr t\int_{\T^p}\dr\mu_t(x)
\inn{\nabla\phi(t,x)}{\frac{1}{h_n}\int_{\R^p}v\g^{h_n}_t(x,v)}\dr v=
\int_a^b\dr t\int_{\T^p}\inn{\nabla\phi(t,x)}{v^{h_n}(t,x)}
\dr\mu_t(x)\tends$$
$$\int_a^b\dr t\int_{\T^p}
\inn{\nabla\phi(t,x)}{X(t,x)}\dr\mu_t(x)  .   \eqno (4.5)$$
As for the Laplacian term, let $i,j\in(1,\dots,p)$; the inequality below comes from H\"older and the fact that 
$\partial^2_{i,j}\phi(t,x)$ is bounded; the limit comes from point 1) of proposition 3.1. 
$$\left\vert
\frac{1}{2h_n}\int_a^b\dr t\int_{\T^p}
\partial^2_{i,j}\phi(t,x)\cdot(h_nv^{h_n}_i(t,x))\cdot
(h_nv^{h_n}_j(t,x))  
\dr\mu_t(x)
\right\vert   \le$$
$$D_2 h_n||v^{h_n}_i||_{L^2(\L^1\otimes\mu_t)}
\cdot
||v^{h_n}_j||_{L^2(\L^1\otimes\mu_t)}\tends 0  .  $$
Together with proposition 3.2, this implies the limit in the formula below, while the second equality comes from the definition of 
$v^{h_n}_i$ and the  third one comes from the definition of the covariance matrix $D^{h_n}=(d^{h_n}_{i,i})_{i,j}$.
$$\frac{1}{2h_n}\int_a^b\dr t\int_{\T^p}\dr\mu_t(x)\int_{\R^p}
\phi^\pprime(t,x)(v,v)\g^{h_n}_t(x,v)\dr v=$$
$$\sum_{i,j}\int_a^b\dr t\int_{\T^p}
\partial^2_{i,j}\phi(t,x)\dr\mu_t(x)\cdot
\frac{1}{2h_n}\int_{\R^p}
v_iv_j\g^{h_n}_t(x,v)\dr v=$$
$$\sum_{i,j}\int_a^b\dr t\int_{\T^p}
\partial^2_{i,j}\phi(t,x)\dr\mu_t(x)\cdot
\frac{1}{2h_n}\int_{\R^p}
(v_i-h_nv^{h_n}_i)(v_j-h_nv^{h_n}_j)\g^{h_n}_t(x,v)\dr v+$$
$$\sum_{i,j}\frac{1}{2h_n}\int_a^b\dr t\int_{\T^p}
\partial^2_{i,j}\phi(t,x)\cdot(h_nv^{h_n}_i)
\cdot(h_nv^{h_n}_j)\dr\mu_t(x)=$$
$$\sum_{i,j}\2\int_a^b\dr t\int_{\T^p}
\partial^2_{i,j}\phi(t,x)\cdot d^{h_n}_{i,j}(t,x)    \dr\mu_t(x)+
\frac{1}{2h_n}\sum_{i,j}\int_a^b\dr t\int_{\T^p}
\partial^2_{i,j}\phi(t,x)\cdot(h_nv^{h_n}_i)\cdot
(h_nv^{h_n}_j)\dr\mu_t(x)
\tends$$
$$\2\int_a^b\dr t\int_{\T^p}\D\phi(t,x)\dr\mu_t(x)  .  $$
Now (4.4) follows from (4.5) and the last formula.

\noindent{\bf Step 2.} Let now ${\cal C}_{(a,b)}(\mu)<+\infty$ and let $\e>0$. By the definition of ${\cal C}_{(a,b)}(\mu)$, we can find a sequence of paths  $\{ \mu_n \}$ converging uniformly to $\mu$ such that, for all $n$, 
$$\Ec_{(a,b)}(\mu_n)\le{\cal C}_{(a,b)}(\mu)+\frac{1}{n}  .  $$
By step 1, there are vector fields $X_n$ such that $\mu_n$ is a weak solution of the Fokker-Planck equation with drift $X_n$; by point 3) of proposition 3.1 and the last formula we have that, for 
$n$ large, 
$$\int_a^b\dr t
\int_{\T^p}\2 |X_n(t,x)|^2\dr\mu_{n,t}(x)\le 
{\cal C}_{(a,b)}(\mu) +\frac{1}{n}  .  \eqno (4.6)$$
Note that $X_n$ and $\mu_n$ induce a one-dimensional current  $T_n$ on $[a,b]\times\T^p$. Indeed, if
$$\o=\o_0\dr t+\o_1\dr x_1+\dots+\o_n\dr x_n=
\o_0\dr t+\o^\prime  $$
is a continuous one-form on $[a,b]\times\T^p$, we can define
$$T_n(\o)=\int_a^b\dr t\int_{\T^p}[
\o_0(t,x)+\o^\prime(t,x)\cdot X_n(t,x)
]   \dr\mu_{n,t}(x)  .  $$
Using (4.6) it is easy to see that the mass norm of $T_n$ is bounded; thus, up to subsequences, $T_n$ converges weakly to a current $T$. In [5] and [10] it is shown how one can define the "kinetic energy" $\phi(S)$ of a current $S$. Actually, they concentrate on closed currents, but the facts we need work even if the current is not closed.  If $S$ is induced by a vector field $X$ and a measure $\mu$ as in the formula above, then $\phi(S)$ has the expression of (4.6), i. e. 
$$\phi(S)=\int_a^b\dr t
\int_{\T^p}\2 |X(t,x)|^2\dr\mu_{t}(x)  .  $$
Now $\phi$ (see again [5] and [10]) is l. s. c. for the weak convergence of currents; thus, (4.6) and the last formula imply that
$$\phi(T)\le{\cal C}_{(a,b)}(\mu)  .  $$
By lemma 3.1 of [10], this implies that $T$ is induced by a vector field $X$ and the measure $\mu_t$; by the last two formulas we have that
$$\int_a^b\dr t
\int_{\T^p}\2 |X(t,x)|^2\dr\mu_{t}\le
{\cal C}_{(a,b)}(\mu)  .  $$
This proves (4.2).

We prove (4.1). Let $\phi\in C^\infty_0((a,b)\times\T^p)$; the first equality below follows from the fact (which we saw at the beginning of this step) that $\mu_n$ is a weak solution of the Fokker-Planck equation with drift $X_n$; the limit comes from the fact that $T_n\tends T$. 
$$0=\int_a^b\dr t\int_{\T^p}[
\D\phi-\inn{\nabla\phi}{X_n}
]     \dr\mu_{n,t}(x)    \tends
\int_a^b\dr t\int_{\T^p}[
\D\phi-\inn{\nabla\phi}{X}
]     \dr\mu_{t}(x)  .  $$
In other words, $\mu_t$ is a solution of the Fokker-Planck equation with drift $X$ and we are done. 

\fin

\vskip 2pc
\centerline{\bf \S 5}
\centerline{\bf Curves which satisfy Fokker-Planck have finite energy}
\vskip 1pc

In this section we are going to end the proof of theorems 1 and 2. We state the converse statement of theorem 1 as a separate proposition.  

\prop{5.1} Let $\mu_t$ be a weak solution of $(FP)_{X}$ and let us suppose that $X\in L^2((a,b)\times\T^p,\L^1\otimes\mu_t)$. Then, 
$${\cal C}_{(a,b)}(\mu_t)\le
\int_a^b\dr t\int_{\T^p}|X(t,x)|^2\dr\mu_t(x)  .   \eqno (5.1)$$

\rm

In lemma 5.2 below, we are going to see that proposition 5.1 holds when the drift is $C^\infty$; the general case will follow by the semicontinuity of lemma 5.3. 

In order to state lemma 5.2 below, we define the cost and forward velocity of a semigroup $P_{s,t}$ on $\Mt$ induced by a Fokker-Planck equation with a sufficiently regular drift. 

\vskip 1pc
\noindent{\bf Definitions.} $\bullet$) Let 
$\g_1\in\dcal_{\mu_1,\mu_2}$ and 
$\g_2\in\dcal_{\mu_2,\mu_3}$; we define $\g_1\ominus\g_2$ as 
$$\g_1\ominus\g_2(x,y)=
\int_{\R^p}\g_2(x+w,y-w)\g_1(x,w)\dr w   .  \eqno (5.2)$$
Let the operation $\ast$ be defined as in section 1; an easy calculation shows that
$$\mu\ast(\g_1\ominus\g_2)=(\mu\ast\g_1)\ast\g_2$$
and that, consequently, 
$\g_1\ominus\g_2\in\dcal_{\mu_1,\mu_3}$. 

\noindent $\bullet$) Let $\{ \g_{\r,\tau} \}_{a\le\r\le\tau\le b}$ be a family of Borel functions on $\T^p\times\R^p$; let $\mu\in\Mt$ and let us set 
$\mu_t=\mu\ast\g_{0,t}$. We say that $\mu_t$ is an orbit of the semigroup $\gam$ starting at $\mu$ if 

\noindent 1) $\g_{s,t}(x,\cdot)\in\dcal_{\mu_s,\mu_t}$ for all 
$a\le s\le t\le b$ and

\noindent 2) if $a\le s_1\le s_2\le s_3\le b$, then
$$\g_{s_1,s_3}=\g_{s_1,s_2}\ominus\g_{s_2,s_3}  .  $$

\noindent $\bullet$) If $\mu_t$ is an orbit of the semigroup $\gam$ starting at $\mu$, we can define as in section 1
$$v_{s,s+h}(x)=\frac{1}{h}\int_{\R^p}v\g_{s,s+h}(x,v)\dr v  ,$$
$$\Ec^h(\mu_s,\mu_{s+h},\{ \g_{\r,\tau} \} )=
\int_{\T^p\times\R^p}A_h(\g_{s,s+h},(x,v))\dr\mu_s(x)\dr v-
\left(\frac{1}{2\pi h}\right)^\frac{p}{2}  ,  $$
$$\Ec^h_{(a,b-h)}(\mu,\{ \g_{\r,\tau} \} ) =
\frac{1}{h}\int_0^{b-h}\Ec^h(\mu_s,\mu_{s+h},\{ \g_{\r,\tau} \} )\dr s  $$
and
$$\Ec_{(a,b)}(\mu,\{ \g_{\r,\tau} \} )=
\liminf_{h\tends 0} \Ec_{(0,b-h)}(\mu,\{ \g_{\r,\tau} \} )  .  $$

\noindent $\bullet$) We can see as in proposition 3.1 that, if 
$\Ec_{(a,b)}(\mu,\{ \g_{\r,\tau} \} )<+\infty$, then there is a sequence $h_n\searrow 0$ and $X\in L^2$ such that
$$v_{s,s+h_n}\weak X\txt{in}
L^2((a,b)\times\T^p,\L^1\otimes\mu_t)  .   $$
We shall say that $X$ is a forward velocity of 
$(\mu,\{ \g_{\r,\tau} \} )$.

\noindent $\bullet$) We say that the orbit $\mu_t$ of the semigroup $\{ \g_{\r,\tau} \}$ solves $(FP)_{X}$ if
the function $\r(t,y)\colon=\g_{s,t}(x,y-x)$ from 
$(s,b)\times\R^p$ to $\R$ is a weak solution of $(FP)_{X}$ for 
$\mu_s$ a. e. $x\in\T^p$ and for all $t>s$. 

\vskip 1pc

We shall need several lemmas.

\lem{5.2} Let $X\in C^\infty([a,b]\times\T^p,\R^p)$, let 
$\mu\in\Mt$ and let $\mu_t$ solve the Fokker-Planck equation with drift $X$ and initial condition $\mu_a=\mu$. Let 
$\{ \g_{\r,\tau} \}$ be the semigroup associated with this equation. Then,

\noindent 1) $X$ is a forward velocity of $(\mu,\{ \g_{\r,\tau} \})$; actually, the vector field $v_{s,s+h}$ defined at the beginning of this section converges to $X$ uniformly.
$${\cal C}_{(a,b)}(\mu_t)\le
\int_a^b\dr t\int_{\T^p}\2 |X(t,x)|^2\dr\mu_t(x)  .   \leqno 2)$$

\proof Point 1) could be proven in a simpler way, but we shall need one of the estimates below for point 2). Let us consider the following stochastic differential equation. 
$$\left\{
\eqalign{
\dr\psi_{s,t}(x)&=X(t,\psi_{s,t}(x))\dr t+\dr w(t)
\txt{for}t\ge s\cr
\psi_{s,s}(x)&=x  .   
}      \right.     $$
Since $\{ \g_{\r,\tau}(x,\cdot) \}$ is the semigroup associated with Fokker-Planck, we easily that $\g_{\r,\tau}(x,\cdot)$ is the law of 
$\psi_{\r,\tau}(x)-x$; together with the definition of push-forward, this implies the last equality below; the expectation ${\bf E}_w$ is for the Wiener measure. The first equality below is the definition of $v_{s,s+h}$.
$$v_{s,s+h}(x)=\frac{1}{h}\int_{\R^p}
v\g_{s,s+h}(x,v)\dr v=\frac{1}{h}{\bf E}_w[
\psi_{s,s+h}(x)-x
]  .   \eqno (5.3)$$
We recall that
$$\psi_{s,s+h}(x)=x+
\int_s^{s+h}X(\tau,\psi_{s,\tau}(x))\dr\tau+w(s+h)-w(s)  .  
\eqno (5.4)$$
On the other side, the Gaussian $N(h X(s,x),hId)$ is the law of 
$\tilde\psi_{s,s+h}(x)-x$, where  $\tilde\psi_{s,s+h}(x)$ satisfies 
$$\tilde\psi_{s,s+h}(x)=x+
\int_s^{s+h}X(s,x)\dr\tau+w(s+h)-w(s)  .  
\eqno (5.5)$$
By well-known properties of the Gaussian, we get the first equality below, while the second one follows by the fact that the Gaussian is the law of $\tilde\psi_{s,s+h}(x)-x$.
$$X(s,x)=\frac{1}{h}
\int_{\R^p}vN(h X(s,x),hId)(v)\dr v=
\frac{1}{h}{\bf E}_w[
\tilde\psi_{s,s+h}-x
]   .  $$
Comparing the last formula and (5.3), we see that point 1) follows if we prove that 
$$\frac{1}{h}{\bf E}_w
|\tilde\psi_{s,s+h}-\psi_{s,s+h}|\tends 0
\txt{as}h\tends 0    \eqno (5.6)$$
uniformly in $s$ and $x$. 
To show this, we subtract (5.5) from (5.4), getting the first inequality below; the third one follows by the fact that $X$ is Lipschitz.
$$|\psi_{s,s+h}(x)-\tilde\psi_{s,s+h}(x)|\le
\int_s^{s+h}|X(\tau,\psi_{s,\tau}(x))-X(s,x)|\dr\tau\le$$
$$\int_s^{s+h}|
X(\tau,\psi_{s,\tau})-X(\tau,\tilde\psi_{s,\tau})
|\dr\tau+
\int_s^{s+h}|
X(\tau,\tilde\psi_{s,\tau})-X(s,x)
|\dr\tau\le  $$
$$D_3\int_s^{s+h}|
\psi_{s,\tau}-\tilde\psi_{s,\tau}
|_{\T^p}\dr\tau+
\int_s^{s+h}|
X(\tau,\tilde\psi_{s,\tau})-X(s,x)
|\dr\tau   .  \eqno (5.7) $$
Let us consider the subset of the Wiener space 
$$A=\{
w\st |\tilde\psi_{s,\tau}(x)-x|>\sqrt\d
\txt{for some}\tau\in[s,s+h]
\}    $$
where we take the distance in $\R^p$, not in $\T^p$.
The first inequality below is Chebyshev, the second one is the standard martingale inequality (see for instance [9], proposition C.5 of the appendix). The equality follows from the fact that the law of $\tilde\psi_{s,t}-x$ is the Gaussian, while the third inequality comes from standard properties of the Gaussian and the fact that $X$ is bounded.
$$\frac{1}{h}{\bf E}_w 1_A\le
\frac{1}{h\d}{\bf E}_w   \left\{
\sup_{\tau\in[0,h]}|\tilde\psi_{s,\tau}-x|^2
\right\}   \le
\frac{4}{h\d}{\bf E}_w|\tilde\psi_{s,h}-x|^2=$$
$$\frac{4}{h\d}\int_{\R^n}
|v|^2 N(X(s,x)h,hId)(v)\dr v\le
\frac{D_4}{\d}   .   $$
The first inequality below comes from (5.7), the fact that $\T^p$ has diameter $\sqrt p$ and the fact that $X$ is bounded; the last one comes from the formula above. 
$$\frac{1}{h}{\bf E}_w[
|\psi_{s,s+h}(x)-\tilde\psi_{s,s+h}(x)|_{\T^p}1_A
]   \le$$
$$\frac{D_3}{h}{\bf E}_w\left[
1_A\int_s^{s+h}\sqrt p\dr\tau
\right]    +
\frac{1}{h}{\bf E}_w\left[
1_A\int_s^{s+h}2||X||_\infty\dr\tau
\right]   \le
D_6{\bf E}_w1_A\le\frac{D_7h}{\d}  .  \eqno (5.8)$$
If $w\not\in A$, (5.7) and the fact that $X$ is Lipschitz imply that, for $\l\in(0,h)$, 
$$|\psi_{s,s+\l}(x)-\tilde\psi_{s,s+\l}(x)|_{\T^p}\le
D_3\int_s^{s+\l}|
\psi_{s,\tau}-\tilde\psi_{s,\tau}
|\dr\tau+D_8h\sqrt\d   .  $$
Since $\psi_{s,s}(x)=\tilde\psi_{s,s}(x)=x$, the Gronwall lemma implies that
$$|\psi_{s,s+h}(x)-\tilde\psi_{s,s+h}(x)|_{\T^p}\le
D_9h\sqrt\d   .  $$
As a consequence, 
$$\frac{1}{h}{\bf E}_w[
|\psi_{s,s+h}(x)-\tilde\psi_{s,s+h}(x)|_{\T^p}\cdot1_{A^c}
]   \le
{\bf E}_w[
D_9\sqrt\d  1_{A^c}
]  \le D_9\sqrt\d   . $$
By the last formula and (5.8), we get that 
$$\frac{1}{h}{\bf E}_w[
|\psi_{s,s+h}(x)-\tilde\psi_{s,s+h}(x)|_{\T^p}
]   \le\frac{D_7h}{\d}+D_9\sqrt\d  .  $$
Given $\e>0$, we can fix $\d$ so small that 
$D_9\sqrt\d<\frac{\e}{2}$; taking $h$ so small that 
$\frac{D_7h}{\d}<\frac{\e}{2}$, formula (5.6) follows. 

We prove point 2). We begin to fix $\e>0$. For 
$s,s+h\in[a+k\e,q+(k+1)\e)$ we define $\tilde\g_{s,s+h}$ as the law of $\tilde\psi_{s,s+h}-x$, where $\tilde\psi_{s,s+h}$ is the solution of (5.5) with drift $X(a+k\e,x)$. We saw above 
that $\tilde\g_{s,s+h}(x,\cdot)$ is a Gaussian. 
For the $a\le s<t\le b$ we define
$$\g^\e_{s,t}=\tilde\g_{s,a+k_1\e}\ominus\tilde\g_{a+k_1\e,a+(k_1+1)\e}
\ominus\dots\ominus\tilde\g_{a+k_l\e,t}$$
where 
$$s\le a+k_1\e<a+(k_1+1)\e<\dots<a+k_l\e\le t   ,   $$
$a+k_1\e$ is the smallest element in $a+\N\e$ larger than $s$, and $a+k_l\e$ the largest one smaller than $t$. It is clear that 
$\g^\e_{s,t}$ defines a semigroup. Said differently, a Dirac delta 
$\d_x$ placed at $x$ at time $a+k\e$ has drift $X(a+k\e,x)$ for 
$t\in[a+k\e,a+(k+1)\e]$. Though the drift is discontinuous in time, it is easy to see that $\mu^\e_t$ is continuous. 

Setting $\mu^\e_t=\mu\ast\g^\e_{s,t}$, it is easy to see (the proof is similar to the one of point 1)) that 
$\mu^\e_t$ converges uniformly to $\mu_t$ as $\e\tends 0$.  Thus, by the definition of ${\cal C}_{a,b}(\mu_t)$, it suffices to show that
$$\lim_{\e\tends 0}\Ec_{(a,b)}(\mu^\e,\{ \g^\e_{\r,\tau} \})\le
\int_a^b\dr t\int_{\T^p}|X(t,x)|^2\dr\mu_t(x)  . \eqno (5.9)$$
To show this, we recall that
$$\Ec^h_{(a,b-h)}(\mu^\e,\{ \g^\e_{\r,\tau} \})=
\frac{1}{h}
\int_a^{b-h}\Ec^h(\mu^\e_s,\mu^\e_{s+h},\{ \g^\e_{\r,\tau} \})\dr s  .  $$
Let us consider $h<\e$ and let us suppose that 
$(s,s+h)\subset[k_i\e,k_{i+1}\e]$; then 
$\g^\e_{s,s+h}(x,v)=N(hX(k_i\e,x),hId)$ (actually, it is a convex combination of Gaussians, but we can forget about this by convexity) and an explicit calculation analogous to the ones in the appendix shows that
$$\Ec^h(\mu^\e_s,\mu^\e_{s+h},\{ \g^\e_{\r,\tau} \})=\frac{h}{2}
\int_{\T^p}
|X(k_i\e,x)|^2\dr\mu^\e_t(x)  .  \eqno (5.10)$$
If $k_i\e\in(s,s+h)$, then $\g^\e_{s,s+h}$ is the convolution of two Gaussians; namely, if $k_i\e=s+h_1$ and $s+h=k_i\e+h_2$, then the Gaussians are $N(h_1X((k_i-1)\e,x),h_1Id)$ and 
$N(h_2X(k_i\e,y),h_2Id)$. Since by our hypotheses $X$ is bounded, another explicit calculation shows that there is $M>0$, independent of $h$ and $\e$, such that
$$\Ec^h(\mu^\e_s,\mu^\e_{s+h},\{ \g^\e_{\r,\tau} \})\le\frac{h}{2}M  .  $$
Since in the interval $(a,b)$ there are at most $\frac{(b-a)}{\e}$ numbers of the form $k_i\e$, the last formula implies that
$$\sum_i 
\int_{k_i\e-h}^{k_i\e}
\frac{1}{h}\Ec^h(\mu^\e_s,\mu^\e_{s+h},\{ \g^\e_{s,t} \})\dr s\le
\frac{h}{2}M\frac{b-a}{\e}  .  $$
We recall that $h<\e$; the first sum below is the contribution of the intervals 
$[s,s+h]$ which do not straddle the points $k_i\e$; the second sum is the contribution of the intervals $[s,s+h]$ straddling some $k_i\e$; the inequality comes from the last formula and (5.10).  
$$\frac{1}{h}
\int_a^{b-h}\Ec^h(\mu^\e_s,\mu^\e_{s+h},\{ \g^\e_{\r,\tau} \})=$$
$$\sum_i\int_{k_i\e}^{(k_i+1)\e-h}
\frac{1}{h}\Ec^h(\mu^\e_s,\mu^\e_{s+h},\{ \g^\e_{\r,\tau} \})\dr s+
\sum_i\int_{k_i\e-h}^{k_i\e}
\frac{1}{h}\Ec^h(\mu^\e_s,\mu^\e_{s+h},\{ \g_{\r,\tau} \})\dr s   \le$$
$$\sum_i\frac{1}{2}\int_{k_i\e}^{(k_i+1)\e -h}\dr t\int_{\T^p}
|X(k_i\e,x)|^2\dr\mu^\e_{t}(x)+
\frac{h}{2}M\frac{b-a}{\e}  .  $$
Letting $h\tends 0$, we get that 
$$\Ec_{(a,b)}(\mu^\e,\{ \g^\e_{\r,\tau} \})\le
\sum_i\int_{k_i\e}^{(k_i+1)\e}\dr t
\int_{\T^p}|X(k_i\e,x)|^2\dr\mu_{t}^\e(x)  .  $$
Letting $\e\tends 0$ and recalling that $X$ is continuous and 
$\mu^\e_t\tends\mu_t$ uniformly, we get (5.9); we saw above that (5.9) implies the thesis.  

\fin

\lem{5.3} The function ${\cal C}_{(a,b)}$ is l. s. c. for uniform convergence. In other words, if $\mu_n\tends\mu$ uniformly on 
$(a,b)$ with respect to the 2-Wasserstein distance, then
$${\cal C}_{(a,b)}(\mu)\le\liminf_{n\tends+\infty}
{\cal C}_{(a,b)}(\mu_n)  .  \eqno (5.11)$$

\proof We recall the stock proof of this fact. By the definition of 
${\cal C}_{(a,b)}(\mu_n)$, we can find curves $\tilde\mu_n$ such that
$$||\mu_n-\tilde\mu_n||_{\sup}\le\frac{1}{n}   \eqno (5.12)$$
and
$$\Ec_{(a,b)}(\tilde\mu_n)\le{\cal C}_{(a,b)}(\mu_n)+\frac{1}{n} .
\eqno (5.13)$$
Since $\mu_n\tends\mu$ uniformly, (5.12) shows that 
$\tilde\mu_n\tends\mu$ uniformly; now the definition of 
${\cal C}_{(a,b)}(\mu)$ implies the first inequality below and (5.13) the second one.
$${\cal C}_{(a,b)}(\mu)\le
\liminf_{n\tends+\infty}\Ec_{(a,b)}(\tilde\mu_n)\le
\liminf_{n\tends+\infty}{\cal C}_{(a,b)}(\mu_n)  .  $$
Since this is (5.11), we are done. 

\fin

\noindent{\bf Proof of proposition 5.1. } Let us define
$$\mu^\e_t=\mu_t\ast N(0,\e Id)  $$
and let us call $\r^\e_t$ the density of $\mu^\e_t$. Here the  Gaussian $N(0,\e Id)$ is in $\R^{p+1}$; since we want a drift of class $C^\infty$ in all variables, we are convoluting also in the time variable. 

Note that, since 
$N(0,\e Id)>0$, also $\r^\e_t$ is strictly positive. Let us consider 
the vector-valued measure
$$ (X_t\mu_t)\ast N(0,\e Id)   .   $$
Since this measure is absolutely continuous with respect to Lebesgue, we call it $E^\e_t\L^p$, with $E^\e_t$ a vector field on 
$\T^p$. We set 
$$X^\e_t=\frac{E^\e_t}{\r^\e_t}  .  $$
Let $\g^\e_{s,t}$ be the semigroup associated with the Fokker-Planck equation of drift $X^\e_t$; let $v^\e_t$ be its forward velocity as defined at the beginning of this section. Since $X^\e\in C^\infty$, we can apply lemma 5.2 and get that 
$$v^\e_t=X^\e_t  .   \eqno (5.14)$$
In step 2 of proposition 4.1 we have defined the functional $\phi$, the "kinetic energy" of a current; since $\phi$ is convex and 
$E_\e$ is a mean of translates of $E$, the inequality below follows from Jensen, while the equalities follow as in step 2 of proposition 4.1.
$$\int_a^b\dr t\int_{\T^p}|X^\e_t(x)|^2\dr\mu^\e_t(x)=
\phi(E_\e)\le\phi(E)=
\int_a^b\dr t\int_{\T^p}|X_t(x)|^2\dr\mu_t(x)   .   \eqno (5.15)$$
The first inequality below is lemma 5.3, the second one follows from the definition of ${\cal C}_{(a,b)}$; the equality is lemma 5.2 while the third inequality comes from the formula above.
$${\cal C}_{(a,b)}(\mu_t)\le
\liminf_{\e\tends 0}{\cal C}_{(a,b)}(\mu^\e_t)\le
\liminf_{\e\tends 0}\Ec_{(a,b)}(\mu^\e_t,\{ \g^\e_{\r,\tau} \})=$$
$$\liminf_{\e\tends 0}
\int_{(a,b)\times\T^p}|X^\e_t(x)|^2\dr\mu^\e_t\dr t\le
\int_{(a,b)\times\T^p}|X_t(x)|^2\dr\mu_t(x)\dr t  .  $$
But this is the thesis.

\fin

\lem{5.4} Let $M>0$, and let $\fun{\mu^n}{[a,b]}{\Mt}$ be curves such that
$${\cal E}_{(a,b)}(\mu^n)\le M\qquad\forall n\ge 1  .
\eqno (5.16)$$
Then, $\{ \mu^n \}$ is compact in $C([a,b],\Mt)$.

\proof We want to use Ascoli-Arzel\`a; since $\Mt$ is compact, it suffices to find a modulus of continuity $\o$ such that
$$d_2(\mu_t,\mu_{t+h})\le\o(h)\qquad
\forall t\in(a,b)    \eqno (5.17)$$
for all curves $\mu$ such that ${\cal E}_{(a,b)}(\mu)\le M$. Let 
$\g_{t,t+h}\in{\cal D}_{\mu_t,\mu_{t+h}}$. We recall from lemma 1.1 that 
$\mu_t\otimes\pi_\sharp(\g_{t,t+h}(x,-x+v)\L^p)$ is a transfer plan between $\mu_t$ and $\mu_{t+h}$; by the definition of $d_2$, this implies the inequality below; the equality follows by the change of variables $\fun{}{y}{v=-x+y}$. 
$$d_2(\mu_t,\mu_{t+h})\le
\left[
\int_{\T^p\times\R^p}|x-y|^2\g_{t,t+h}(x,-x+y)\dr\mu_t(x)\dr y
\right]^\2     =$$
$$\left[
\int_{\T^p\times\R^p}|v|^2\g_{t,t+h}(x,v)\dr\mu_t(x)\dr v 
\right]^\2 .  $$
Thus, (5.17) follows if we prove that, if $(\mu,\{ \g_{s,t} \})$ is a weak solution of the Fokker-Planck equation with drift $X$, then
$$\int_{\T^p\times\R^p}|v|^2\g_{t,t+h}(x,v)\dr\mu_t(x)\dr v\le
2 h\int_a^b\dr\tau\int_{\T^p}|X(\tau,x)|^2\dr\mu_\tau(x)+2 h  .  
\eqno (5.18)  $$
To prove this, let us begin to suppose that 
$X\in C^\infty([a,b]\times\T^p)$; then for $t\ge s$ the measure 
$\mu_t$ is the law of the solution $x(t)$ of
$$\dr x(t)=X(t,x(t))\dr t+\dr w(t)  \eqno (5.19)$$
where the initial condition $x(s)$ has law $\mu_s$ and $w(t)$ is a Brownian motion on $(a,+\infty)$. In other words,
$$x(t)=x(s)+
\int_s^tX(\tau,x(\tau))\dr\tau+w(t)-w(s) $$
where $x(s)$ has law $\mu_s$ and is independent from $\mu_t$ for $t>s$. Let us denote as usual by ${\bf E}_w$ the expectation with respect to the Wiener measure. The first equality below comes from the fact that $\g_{s,t}$ is the semigroup induced by (5.19), the first inequality comes by the formula above and 
H\"older, the second one is H\"older and the last equality comes from well-known properties of the Brownian motion.
$$\int_{\T^p\times\R^p}|v|^2\g_{s,t}(x,v)\dr\mu_s(x)\dr v=
{\bf E}_w |x(t,w)-x(s,w)|^2\le$$
$$2{\bf E}_w
\left\vert
\int_s^t X(\tau,x(\tau))\dr\tau
\right\vert^2+
2{\bf E}_w
|w(t)-w(s)|^2\le$$
$$2(t-s){\bf E}_w
\int_s^t |X(\tau,x(\tau))|^2\dr\tau+
2{\bf E}_w
|w(t)-w(s)|^2=$$
$$2(t-s)\int_s^t\dr\tau\int_{\T^p}|X(\tau,x)|^2\dr\mu_\tau(x)+
2(t-s).  $$
But this is (5.18) for the smooth drift $X$. 

We prove the general case. Let us approximate $X$ with smooth vector fields $X_\e$ as in the proof of proposition 5.1; let us call 
$\{ \g^\e_{s,t} \}$ the semigroup of the Fokker-Planck equation with drift $X_\e$. The first inequality below follows from the lower semicontinuity of the functional
$$\fun{}{\nu}{\int_{\T^p\times\R^p}|v|^2\dr\nu(x,v)}$$
under weak convergence. The second one is (5.18) for the smooth drift $X_\e$ and the third one follows as in (5.15). 
$$\int_{\T^p\times\R^p}|v|^2\g_{t,t+h}(x,v)\dr\mu_t(x)\dr v\le
\lim_{\e\tends 0}\int_{\T^p\times\R^p}
|v|^2\g^\e_{t,t+h}(x,v)\dr\mu^\e_t(x)\dr v  \le  $$
$$\lim_{\e\tends 0}\left[
2 h\int_a^b\dr\tau\int_{\T^p}|X_\e(\tau,x)|^2\dr\mu^\e_\tau(x)+
2 h
\right]  \le  
2 h\int_a^b\dr\tau\int_{\T^p}|X(\tau,x)|^2\dr\mu_\tau(x)+
2 h   .  $$
Since this is (5.18), for the drift $X$, we are done.

\fin

\noindent{\bf Remark.} By the last lemma, if  
$\Ec_{(a,b)}(\mu)<+\infty$, then the curve of measures $\mu$ is continuous. This allows us to embed the initial condition in the definition of weak solution of Fokker-Planck; namely, $\mu$ satisfies 
$$\int_a^b\dr t\int_{\T^p}
\partial_t\phi(t,x)\dr\mu_t(x)+\int_{\T^p}\phi(a,x)\dr\mu_a(x)=$$
$$-\2\int_a^b\dr t\int_{\T^p} \D\phi(t,x)\dr\mu_t(x)-
\int_a^b\dr t\int_{\T^p}\inn{\nabla\phi(t,x)}{X(t,x)}\dr\mu_t(x)
\quad\forall\phi\in C^\infty_c([a,b)\times\T^p) .  $$
We omit the proof of this, since it follows in a standard way from (4.1) and the continuity of $\mu$. 

\vskip 1pc
\noindent{\bf End of the proof of theorem 1.} Let $\mu_t$ be a curve of measures as in the hypotheses of theorem 1. By proposition 4.1 $\mu$ solves $(FP)_{X}$, while (4.2) proves half of equality (3). The converse, and the opposite inequality of (4.2), follows by proposition 5.1.

\vskip 1pc
\noindent{\bf End of the proof of theorem 2.} Point 1) is lemma 5.4, while point 2) is lemma 5.3.

\vskip 2pc
\centerline{\bf Appendix}
\centerline{\bf Estimates on the Gaussian}
\vskip 1pc

In this appendix, we prove the estimates on the Gaussian we use  throughout the paper. For starters, we fix $h>0$ and give some definitions.

\vskip 1pc 
\noindent{\bf Definitions.} $\bullet$) First of all, we settle the notation for the Gaussian: if $Q$ is a symmetric, positive-definite matrix and 
$a\in\R^p$, we denote the Gaussian of mean $a$ and variance 
$Q$ by
$$N(a,Q)(v)=\left(
\frac{1}{(2\pi)^p\det Q}
\right)^\2
e^{-\2\inn{Q^{-1}(v-a)}{v-a}}  .   $$

\noindent $\bullet$) In section 1, we have defined $Den$ as the set of all Borel probability densities on $\R^p$; here, we further define  $Den_2$ as the set of all the Borel probability densities on $\R^p$ whose second moments are finite.

\noindent $\bullet$) Let $a\in\R^p$; we group in a set $Mean(a)$ the functions $\g\in Den_2$ such that
$$\frac{1}{h}\int_{\R^p}v\g(v)\dr v=a  .  $$
Note that the integral converges, i. e. $v\in L^1(\g\L^p)$: this follows by H\"older's inequality since $1\in L^1(\g\L^p)$ ($\g$ is a probability density) and 
$|v|^2\in L^1(\g\L^p)$ ($\g$ has finite second moment).

\noindent $\bullet$) For $i,j\in(1,\dots,p)$, $a\in\R^p$ and $\d>0$, we define $Corr_{i,j}(a,\d)$ as the set of the functions 
$\g\in Mean(a)$ such that
$$\frac{1}{h}\int_{\R^p}
(v_i-ha_i)(v_j-ha_j)\g(v)\dr v=\d   $$

\noindent $\bullet$) For $a\in\R^p$ and $\d>0$, we define 
$Trace(a,\d)$ as the set of the $\g\in Mean(a)$ such that
$$\frac{1}{h}\int_{\R^p}|v-ah|^2\g(v)\dr v=p\d  .  $$

\noindent $\bullet$) For $a\in\R^p$ and $\d>0$, we define
$$B_{off-diag}(a,\d)=
\frac{-1+\sqrt{1+4\d^2}}{2}+\frac{h}{2}|a|^2+
\log\left[
\frac{-1+\sqrt{1+4\d^2}}{2\d^2(2\pi h)^p}
\right]^\2   ,  $$
$$B_{diag}(a,\d)=\frac{\d-1}{2}+\frac{h}{2}|a|^2+
\log\left(
\frac{1}{2\pi h\d}
\right)^\2+
\log\left(
\frac{1}{2\pi h}
\right)^\frac{p-1}{2}      $$
and
$$T(a,\d)=p\frac{\d-1}{2}+\frac{h}{2}|a|^2+
\log\left(
\frac{1}{2\pi h\d}
\right)^\frac{p}{2}    .  $$
Our first lemma is an estimate on the trace.

\lem{A.1} Let $a\in\R^p$, $\d>0$ and let $\g\in Trace(a,\d)$. Let 
$A_h$ be defined as in section 1. Then, 
$$\int_{\R^p}A_h(\g,v)\dr v\ge T(a,\d)   .   \eqno (A.1)$$
Moreover, if $\g\in Den\setminus Den_2$, then
$$\int_{\R^p}A_h(\g,v)\dr v=+\infty  .  \eqno (A.2)$$

\proof We begin with (A.1).  Our plan is to consider the functional 
$$\fun{I}{\g}{\int_{\R^p}A_h(\g,v)\dr v}$$
and minimize it over $Trace(a,\d)$; we shall show that the minimal $\g$ exists and that $I(\g)=T(a,\d)$. 

We begin to note that, since $I$ is strictly convex and the set 
$Trace(a,\d)$ is convex, there is at most one minimizer. It is standard ([13] or proposition I, 5.6 of [11]) that, if we find 
$\g\in L^1((1+\2|v|^2)\L^p)$, $b\in\R^p$ and $\eta,\l\in\R$ which solve the Lagrange multiplier problem
$$\left\{
\eqalign{
\cinh{v}+1+\log\g(v)&=\l+\inn{b}{v}+\frac{\eta}{2h}|v-ah|^2\cr
\int_{\R^p}\g(v)\dr v&=1\cr
\frac{1}{h}\int_{\R^p}v\g(v)\dr v&=a\cr
\frac{1}{h}\int_{\R^p}|v-ah|^2\g(v)\dr v&=p\d
}     \right.    \eqno (A.3)$$
then $\g$ is the unique minimizer of $I$ on $Trace(a,\d)$. From the first equation of (A.3) we get the first equality below.
$$\g(v)=e^{\l-1}\exp\left(
-\cinh{v}+\frac{\eta}{2h}|v-ah|^2+\inn{b}{v}
\right)   =   $$
$$\exp\left(
\l-1+\frac{\eta h}{2}|a|^2+\frac{h}{2(1-\eta)}|b-\eta a|^2
\right)   \cdot
\exp\left(
-\frac{1-\eta}{2h}
\left\vert
v-\frac{h}{1-\eta}(b-\eta a)
\right\vert^2
\right)    .  $$
This is a Gaussian multiplied by a complicated coefficient; the second equation of (A.3) makes short work of it.
$$\g(v)=\left(
\frac{1-\eta}{2\pi h}
\right)^\frac{p}{2}
\exp\left(
 -\frac{1-\eta}{2h}
\left\vert
v-\frac{h}{1-\eta}(b-\eta a)
\right\vert^2
\right)   .   $$
Together with the third formula of (A.3), this implies that
$$\frac{h}{1-\eta}(b-\eta a)=ah  .  $$
Substituting this into the expression for $\g$, we get that
$$\g(v)=\left(
\frac{1-\eta}{2\pi h}
\right)^\frac{p}{2}
\exp\left( 
-\frac{1-\eta}{2h}
\left\vert
v-ah
\right\vert^2 
\right)  .   $$
Together with the fourth formula of (A.3), this implies that
$$p\frac{h}{1-\eta}=ph\d    $$
which plugged into the expression for $\g$  yields
$$\g(v)=\left(
\frac{1}{2\pi h\d}
\right)^\frac{p}{2}
\exp\left(
-\frac{1}{2h\d}
|v-ah|^2
\right)   .   $$
Substituting, we get the fourth equality below, while the third and the last one come from (A.3).
$$I(\g)=\int_{\R^p}\left[
\cinh{v}+\log\g(v)
\right]  \g(v)\dr v=$$
$$\int_{\R^p}\left[
\cinh{v-ah}+\inn{a}{v}-\frac{h}{2}|a|^2+
\log\g(v)
\right]   \g(v)\dr v  =
\frac{p\d}{2}+\frac{h}{2}|a|^2+\int_{\R^p}\g(v)\log\g(v)\dr v=$$
$$\frac{p\d}{2}+\frac{h}{2}|a|^2+\int_{\R^p}
\left[
\log\left(
\frac{1}{2\pi h\d}
\right)^\frac{p}{2}-
\frac{1}{2h\d}|v-ah|^2
\right]   \g(v)\dr v=$$
$$\frac{p\d}{2}+\frac{h}{2}|a|^2+
\log\left(
\frac{1}{2\pi h\d}
\right)^\frac{p}{2}-\frac{p}{2}  .   $$
Since this is (A.1), we are done. 

We prove (A.2). Let $\g\in Den\setminus Den_2$; since the function
$$\fun{}{t}{
\left(
\cinh{v}+\log t
\right)   t
}  $$
is convex, its graph is above its tangent at 
$t=e^{-\frac{1}{4h}|v|^2}$, i. e. the inequality below holds. 
$$\left[
\cinh{v}+\log t
\right] t\ge
\frac{1}{4h}|v|^2e^{
-\frac{1}{4h}|v|^2
}      +
\left[
\frac{1}{4h}|v|^2+1
\right]\cdot
\left[
t-e^{
-\frac{1}{4h}|v|^2
} 
\right]   =$$
$$-e^{
-\frac{1}{4h}|v|^2
} +
\left[
\frac{1}{4h}|v|^2+1
\right]\cdot   t   .  $$
This implies the inequality below.
$$\int_{\R^p}\left[
\cinh{v}+\log\g(v)
\right]  \g(v)\dr v\ge$$
$$\int_{\R^p}
\left[
-e^{
-\frac{1}{4h}|v|^2
}    +\left(
\frac{1}{4h}|v|^2+1
\right)   \g(v)
\right]     \dr v   .  $$
Let us suppose by contradiction that $I(\g)<+\infty$; if in the formula above we get rid of the terms which are obviously finite, we get that  
$$\int_{\R^p}
\frac{1}{4h}|v|^2\g(v)\dr v<+\infty    $$
contradicting the fact that 
$\g\in Den\setminus Den_2$.

\fin

We also need an estimate on each single element of the covariance matrix.

\lem{A.2} Let $i,j\in(1,\dots,p)$. Then,
$$\min\left\{
\int_{\R^p}A_h(\g,v)\dr v\st \g\in Corr_{i,j}(a,\d)
\right\}   \ge
\left\{
\eqalign{
B_{diag}(a,\d)&\txt{if} i=j\cr
B_{off-diag}(a,\d)&\txt{if}i\not =j   .  
}       \right.    \eqno(A.4)   $$

\proof We begin with the case in which $i=j$; without loss of generality we can suppose that $i=j=1$. As in lemma A.1, we are going to write down explicitly the minimal for the left hand side of (A.4); we look for $\g\in L^1((1+\2|v|^2)\L^p)$ and Lagrange multipliers $b\in\R^p$, $\l ,\eta\in\R$ which solve
$$\left\{
\eqalign{
\cinh{v}+1+\log\g(v)&=\l+\frac{\eta}{2h}|v_1-a_1h|^2+\inn{b}{v}\cr
\int_{\R^p}\g(v)\dr v&=1\cr
\frac{1}{h}\int_{\R^p}v\g(v)\dr v&=a\cr
\frac{1}{h}\int_{\R^p}|v_1-a_1h|^2\g(v)&=\d  .
}     \right.        \eqno (A.5)$$
From the first equation of (A.5) we get the first equality below; for the second one, we have set $v^\prime=(v_2,\dots,v_p)$.
$$\g(v)=e^{\l-1}\exp
\left(
-\cinh{v}+\frac{\eta}{2h}|v_1-a_1h|^2+\inn{b}{v}
\right)=$$
$$e^{\l-1}\exp
\left(
-\cinh{v_1}+\frac{\eta}{2h}|v_1-a_1h|^2+b_1v_1
\right)   \cdot
\exp\left(
-\cinh{v^\prime}+\inn{b^\prime}{v^\prime}
\right)  =  $$
$$\exp\left(
\l-1+\frac{\eta h}{2}|a_1|^2+\frac{h}{2(1-\eta)}|b_1-\eta a_1|^2
+\frac{h}{2}|b^\prime|^2
\right)   \cdot$$
$$\exp\left(
-\frac{1-\eta}{2h}\left\vert
v_1-\frac{h}{1-\eta}(b_1-\eta a_1)
\right\vert^2
\right)   \cdot
\exp\left(
-\cinh{v^\prime-hb^\prime}
\right)   .  $$
This is a Gaussian multiplied by a complicated expression; the second formula of (A.5) makes short work of it.
$$\g(v)=\left(
\frac{1-\eta}{2\pi h}
\right)^\2\cdot
\left(
\frac{1}{2\pi h}
\right)^\frac{p-1}{2}\cdot
\exp\left(
-\frac{1-\eta}{2h}\left\vert
v_1-\frac{h}{1-\eta}(b_1-\eta a_1)
\right\vert^2-\cinh{v^\prime-hb^\prime}
\right)   .  $$
As in the last lemma, this formula and the third one of (A.5) give us two different expressions for the mean of $\g$:
$$\left(
\frac{h}{1-\eta}(b_1-\eta a_1),hb^\prime
\right)=(ha_1,ha^\prime)   .   $$
Substituting into the expression for $\g$, we get that
$$\g(v)=\left(
\frac{1-\eta}{2\pi h}
\right)^\2\cdot
\left(
\frac{1}{2\pi h}
\right)^\frac{p-1}{2}\cdot
\exp\left(
-\frac{1-\eta}{2h}\left\vert
v_1-a_1h
\right\vert^2-\cinh{v^\prime-ha^\prime}
\right)   .  $$
By this formula and the fourth one of (A.5), we can write in two different ways the variance of $\g$ in the $x_1$ direction:
$$\frac{h}{1-\eta}=\d h  .   $$
Substituting into the expression for $\g$, we get that
$$\g(v)=\left(
\frac{1}{2\pi h\d}
\right)^\2\cdot
\left(
\frac{1}{2\pi h}
\right)^\frac{p-1}{2}\cdot
\exp\left(
-\frac{1}{2h\d}\left\vert
v_1-a_1h
\right\vert^2-\cinh{v^\prime-ha^\prime}
\right)   .  $$
We use this to get the fourth equality below; the third and the last one come from (A.5).
$$I(\g)=\int_{\R^p}\left[
\cinh{v}+\log\g(v)
\right]   \g(v)\dr v=$$
$$\int_{\R^p}\left[
\cinh{v_1-a_1h}+\cinh{v^\prime-a^\prime h}+
a_1v_1+\inn{a^\prime}{v^\prime}-\frac{h}{2}|a_1|^2
-\frac{h}{2}|a^\prime|^2+\log\g(v)
\right]     \g(v)\dr v=$$
$$\frac{\d}{2}+\frac{p-1}{2}+\frac{h}{2}|a|^2+
\int_{\R^p}\g(v)\log\g(v)\dr v=$$
$$\frac{\d}{2}+\frac{p-1}{2}+\frac{h}{2}|a|^2+
\int_{\R^p}\left[
-\frac{1}{2h\d}|v_1-a_1h|^2-\cinh{v^\prime-a^\prime h}
+\log\left(
\frac{1}{2\pi h\d}
\right)^\2+
\log\left(
\frac{1}{2\pi h}
\right)^\frac{p-1}{2}
\right]    \g(v)\dr v=$$
$$\frac{\d-1}{2}+\frac{h}{2}|a|^2
+\log\left(
\frac{1}{2\pi h\d}
\right)^\2+
\log\left(
\frac{1}{2\pi h}
\right)^\frac{p-1}{2}    .   $$
By the definition of $B_{diag}$, this is the first inequality of (A.4).

Now we tackle the off-diagonal case, $i\not=j$. Without loss of generality, we can suppose that $i=1$ and $j=2$. As before, we have to find $\g\in L^1((1+\cin{v})\L^p)$ and Lagrange multipliers $\l,\eta\in\R$, $b\in\R^p$ such that
$$\left\{
\eqalign{
\cinh{v}+1+\log\g(v)&=\l+\frac{\eta}{h}(v_1-a_1h)(v_2-a_2h)+
\inn{b}{v}\cr
\int_{\R^p}\g(v)\dr v&=1\cr
\frac{1}{h}\int_{\R^p}v\g(v)\dr v&=a\cr
\frac{1}{h}\int_{\R^p}(v_1-a_1h)(v_2-a_2h)\g(v)\dr v&=\d  . 
}
\right.   \eqno (A.6)$$
From the first equation of (A.6) we get the first equality below.
$$\g(v)=e^{\l-1}\exp\left(
-\cinh{v}+\frac{\eta}{h}(v_1-a_1h)(v_2-a_2h)+\inn{b}{v}
\right)=$$
$$\exp\left(
\l-1+\frac{h}{2}|a|^2
\right)   \cdot
\exp\left(
-\cinh{v-ah}+\frac{\eta}{h}(v_1-a_1h)(v_2-a_2h)+
\inn{b-a}{v}
\right)   .   $$
Now we set 
$$Q^{-1}=Id-\eta(e_1\otimes e_2+e_2\otimes e_1)    \eqno (A.7)$$
where $Id$ is the identity matrix and $\{ e_i \}$ is the standard basis of $\R^p$; the last formula becomes 
$$\g(v)=\exp\left(
\l-1+\frac{h}{2}|a|^2+\inn{ah}{b-a}
\right)   \cdot
\exp\left(
-\frac{1}{2h}\inn{Q^{-1}(v-ah)}{v-ah}+\inn{b-a}{v-ah}
\right)   .   $$
Naturally, at the end we shall have to check that $\eta\in(-1,1)$, since $Q$ must be positive-definite if we want $\g$ to be integrable. Setting
$$\bar b=Q(b-a)$$
we can write
$$\g(v)=\exp\left(
\l-1+\frac{h}{2}|a|^2+\inn{ah}{b-a}+
\frac{h}{2}\inn{Q^{-1}\bar b}{\bar b}
\right)   \cdot
\exp\left(
-\frac{1}{2h}\inn{Q^{-1}(v-(a+\bar b)h)}{v-(a+\bar b)h}
\right)   .   $$
The second formula of (A.6) settles the constant before the Gaussian:
$$\g(v)=\left(
\frac{1}{(2\pi h)^p\det Q}
\right)^\2
\exp\left(
-\frac{1}{2h}\inn{Q^{-1}(v-(a+\bar b)h)}{v-(a+\bar b)h}
\right)   .   $$
Since (A.7) implies that
$${\rm det}Q^{-1}=1-\eta^2$$
the last formula becomes
$$\g(v)=\left(
\frac{1-\eta^2}{(2\pi h)^p}
\right)^\2
\exp\left(
-\frac{1}{2h}\inn{Q^{-1}(v-(a+\bar b)h)}{v-(a+\bar b)h}
\right)   .   $$
As in the first part of this lemma, this and the third formula of (A.6) give two different expressions for the mean of $\g$:
$$(a+\bar b)h=ah   .   $$
Thus,
$$\g(v)=\left(
\frac{1-\eta^2}{(2\pi h)^p}
\right)^\2
\exp\left(
-\frac{1}{2h}\inn{Q^{-1}(v-ah)}{v-ah}
\right)   .   \eqno (A.8)$$
This is a Gaussian whose covariance matrix is $hQ$; the last one of (A.6) implies that $Q_{1,2}=\d$. Inverting (A.7) is another way of calculating $Q_{1,2}$; a quick calculation now yields
$$\frac{\eta}{1-\eta^2}=\d  $$
or, equivalently,
$$\eta=\frac{-1+\sqrt{1+4\d^2}}{2\d},\qquad
\a\colon=1-\eta^2=\frac{-1+\sqrt{1+4\d^2}}{2\d^2}   .   
\eqno (A.9)$$
Solving the second degree equation for $\eta$ we have chosen the positive square root because with this choice $\eta\in(-1,1)$ and $Q$ is positive-definite, as we wanted.

With $\a$ defined as in (A.9), (A.8) becomes
$$\g(v)=\left[
\frac{\a}{(2\pi h)^p}
\right]^\2
\exp\left(
-\frac{1}{2h}\inn{Q^{-1}(v-ah)}{v-ah}
\right)   .   $$
Together with the third formula of (A.6), this yields the third equality below; the fourth one comes from the definition of $Q^{-1}$ in (A.7) and the last equality comes from (A.6).
$$I(\g)=\int_{\R^p}\left[
\cinh{v}+\log\g(v)
\right]   \g(v)\dr v=$$
$$\int_{\R^p}\left[
\cinh{v-ah}+\inn{a}{v}-\frac{h}{2}|a|^2+\log\g(v)
\right]   \g(v)\dr v=$$
$$\frac{h}{2}|a|^2+\int_{\R^p}\left\{
\cinh{v-ah}+\log\left[
\frac{\a}{(2\pi h)^p}
\right]^\2-
\frac{1}{2h}\inn{Q^{-1}(v-ah)}{v-ah}
\right\}   \g(v)\dr v=$$
$$\frac{h}{2}|a|^2+\int_{\R^p}\left\{
\frac{\eta}{h}(v_1-a_2h)(v_2-a_2h)+
\log\left[
\frac{\a}{(2\pi h)^p}
\right]^\2
\right\}   \g(v)\dr v=$$
$$\frac{h}{2}|a|^2+\d\eta+\log\left[
\frac{\a}{(2\pi h)^p}
\right]^\2   .   $$
By (A.9) and the definition of $B_{off-diag}$ at the beginning of this section, this implies the second inequality of (A.4).

\fin

\lem{A.3} There is $D_1>0$, independent of $a\in\R^p$ and $h>0$, such that the following holds.

\noindent 1) The function $\fun{}{\d}{T(a,\d)}$ from $(0,+\infty)$ to 
$\R$ reaches its minimum at $\d=1$; moreover,
$$T(a,\d)-T(a,1)\ge\left\{
\eqalign{
D_1|\d-1|^2&\txt{if}\d\le 2\cr
D_1(\d-1)&\txt{if}\d>2   .    
}     \right.      \eqno (A.10)$$

\noindent 2) The function $\fun{}{\d}{B_{diag}(a,\d)}$ from 
$(0,+\infty)$ to $\R$ reaches its minimum for $\d=1$ and 
$$B_{diag}(a,\d)-B_{diag}(a,1)\ge\left\{
\eqalign{
D_1|\d-1|^2&\txt{if}\d\le 2\cr
D_1(\d-1)&\txt{if}\d>2   .    
}     \right.      \eqno (A.11)$$

\noindent 3) The function $\fun{}{\d}{B_{off-diag}(a,\d)}$ from 
$[0,+\infty)$ to $\R$ reaches its minimum for $\d=0$; moreover,
$$B_{off-diag}(a,\d)-B_{off-diag}(a,0)\ge\left\{
\eqalign{
D_1\d^2&\txt{if}\d\le 1\cr
D_1\d&\txt{if}\d>1   .    
}     \right.      \eqno (A.12)$$

\proof Since the proof of (A.10) is analogous to that of (A.11), we prove the former. We set
$$g(\d)=p\frac{\d-1}{2}-\frac{p}{2}\log\d   $$
and note that, by the definition of $T(a,\d)$, 
$$T(a,\d)-T(a,1)=g(\d)   .  $$
Thus, it suffices to show that the minimum of $g$ is in $\d=1$, that 
$g(1)=0$ and that $g$ satisfies (A.10). These assertions follow from freshman analysis; we prove the last one. Since
$$g^\pprime(\d)=\frac{p}{2\d^2}  ,  $$
we get that $\fun{g}{(0,+\infty)}{\R}$ is strictly convex; we also get that $g^\pprime(\d)\ge\frac{p}{8}$ if $\d\in(0,2)$, so that the first inequality of (A.10) holds; the second one follows recalling that, since $g$ is convex, its derivative in $[2,+\infty)$  is larger than 
$g^\prime(2)$, which is positive.

We prove (A.12). We begin to set 
$$l(\d^2)=B_{off-diag}(a,\d)-B_{off-diag}(a,0)$$
so that
$$l(t)=\frac{-1+\sqrt{1+4t}}{2}+
\2\log(\sqrt{1+4t}-1)-\2\log(t)-\2\log 2  .  $$
Instead of studying $l$ for $t\in(0,+\infty)$, we set 
$$s=\sqrt{1+4t}$$
and study
$$m(s)\colon=l(\frac{s^2-1}{4})=\frac{s-1}{2}+\2\log(s-1)-
\2\log(\frac{s^2-1}{4})-\2\log 2=
\frac{s-1}{2}-\2\log(s+1)+\2\log 2   $$
for $s\in(1,+\infty)$. We get that
$$m^\prime(s)=\frac{s}{2(s+1)}$$
which shows that $m$ is monotone increasing in $[1,+\infty)$; going back to the variable $t$, we get that $l$ is monotone increasing on 
$[0,+\infty)$. By the formula above, $m$ is convex and 
$m^\prime(1)=\frac{1}{4}$; since $m(1)=0$, this implies that 
$$m(s)\ge\frac{1}{4} (s-1)   \txt{for}  s\ge 1   .   $$
Recalling that $s=\sqrt{1+4t}$, we get that
$$l(t)\ge\frac{1}{4} (\sqrt{1+4t}-1)  \txt{for}t\ge 0  $$
i. e. that
$$l(\d^2)\ge\frac{1}{4} (\sqrt{1+4\d^2}-1)  \txt{for}\d\ge 0  $$
which immediately implies (A.12).

\fin

We need a lemma to estimate the contribution of $\g|_{B(0,r)^c}$ to the second moment; we begin with two definitions.

\vskip 1pc

\noindent {\bf Definitions.} $\bullet$) For $r,\d>0$ we define 
$Out(r,\d)$ as the subset of the $\g\in Den$ such that
$$\frac{1}{2h}\int_{B(0,r)^c}|v|^2\g(v)\dr v=\d  .  $$

\noindent $\bullet$ We define
$$F_h(r,\d)=\inf\left\{
\int_{\R^p}A_h(\g,v)\dr v-
\log\left(
\frac{1}{2\pi h}
\right)^\frac{p}{2}\st  \g\in Out(r,\d)
\right\}   .   \eqno (A.13)$$

\lem{A.4} 1) Let $r>0$; then there is $\d_0(r,h)\tends 0$ as 
$h\tends 0$ such that 
$$F_h(r,\d)\ge
\frac{\d}{2}  \txt{if} \d\ge\d_0(r,h)  .  $$

\proof As in lemmas A.1 and A.2, the function $\g$ which minimizes in (A.13) is the solution of the following Lagrange multiplier problem.
$$\left\{
\eqalign{
\cinh{v}+1+\log\g(v)&=\l-\frac{\eta}{2h}|v|^2 1_{B(0,r)^c}(v)\cr
\int_{\R^p}\g(v)\dr v&=1\cr
\frac{1}{2h}\int_{B(0,r)^c}|v|^2\g(v)\dr v&=\d
}    \right.     \eqno (A.14)  $$
From the first equation of (A.14) we get that
$$\g_{\l,\eta}=e^{\l-1}\exp\left(
-\cinh{v}-\frac{\eta}{2h}|v|^21_{B(0,r)^c}(v)
\right)   .    \eqno (A.15)$$

\noindent{\bf Step 1.} We assert that there are Lagrange multipliers $(\l,\eta)$ such that the function $\g_{\l,\eta}$ of (A.15) satisfies the second and third equations of (A.14). Before proving existence, we note that $(\l,\eta)$ will be unique, because of the strict convexity of the functional. 

As for existence, let us note that for all $\eta>-1$ there is 
$\l(\eta)\in\R$ for which the second equation of (A.14) holds: by (A.15) it suffices to take
$$e^{1-\l(\eta)}=\int_{\R^p}
\exp\left(
-\cinh{v}-\frac{\eta}{2h}|v|^21_{B(0,r)^c}(v)
\right) \dr v   .   \eqno (A.16)$$
It remains to show that we can find $\eta(\d)$ such that 
$\g_{\l(\eta(\d)),\eta(\d)}$ satisfies the third equation of (A.14). To prove this, we define the function 
$$\fun{\bar\d}{\eta}{
\frac{1}{2h}\int_{B(0,r)^c}|v|^2\g_{\l(\eta),\eta}(v)\dr v
}    .  $$
This is clearly a continuous function of $\eta$; we want to use the intermediate value theorem to show that the equation 
$\bar\d(\eta)=\d$ has a solution $\eta$. We begin to see that dominated convergence and (A.16) yield the first equality below, while the second one is the definition of $a(h)$. 
$$\lim_{\eta\tends+\infty} e^{
1-\l(\eta)
}=  
\int_{B(0,r)}e^{-\cinh{v}}\dr v
\colon =   (2\pi h)^\frac{p}{2}[1-a(h)]  .  \eqno (A.17)$$
Note that, by well-known properties of the Gaussian, 
$$0\le a(h)\le e^{-\frac{D_2}{ h}}  \txt{for}h\in(0,1]  .  
\eqno (A.18)$$
Here and in the following we denote by $D_i$ a constant depending only on $r>0$ and we forget the dependence of $a(h)$ on $r$; we can do this because $r$ is fixed throughout the lemma. 

The first equality below follows from the definition of the function 
$\bar\d$ and (A.15), the second one comes from (A.17) and the  last one is dominated convergence. 
$$\lim_{\eta\tends+\infty}\bar\d(\eta)=
\lim_{\eta\tends+\infty}\frac{1}{2h}
\int_{B(0,r)^c}|v|^2e^{\l(\eta)-1}
e^{-\frac{1+\eta}{2h}|v|^2}   \dr v=$$
$$\lim_{\eta\tends+\infty}\frac{1}{2h}\int_{B(0,r)^c}
|v|^2
\frac{1}{(2\pi h)^\frac{p}{h}[1-a(h)]}
e^{
-\frac{1+\eta}{2h}|v|^2
}    \dr v   =0  .  $$
On the other side, the definition of $\bar\d$ implies the first equality below, while (A.15) and (A.16) imply the second one; the third one follows by the change of variables 
$y=\sqrt\frac{1+\eta}{h}v$. As for the last equality, note that the three integrals in the expression on the left tend to a positive limit as $\eta\searrow -1$, but $\left(
\frac{h}{1+\eta}
\right)^{\frac{p}{2}+1}$ goes to $+\infty$ faster that $\left(
\frac{h}{1+\eta}
\right)^{\frac{p}{2}}$. 
$$\lim_{\eta\searrow -1}\bar\d(\eta)=
\lim_{\eta\searrow -1}\frac{1}{2h}
\int_{B(0,r)^c}|v|^2\g_{\l(\eta),\eta}(v)\dr v=$$
$$\lim_{\eta\searrow -1}\frac{1}{2h}\cdot
\frac{
\int_{B(0,r)^c}|v|^2e^{-\frac{1+\eta}{2h}|v|^2}\dr v
}{
\int_{\R^p}\exp\left[
-\cinh{v}-\frac{\eta}{2h}|v|^21_{B(0,r)^c}(v)
\right]  \dr v
}  =  $$
$$\lim_{\eta\searrow -1}\frac{1}{2h}
\cdot
\frac{   
\left(
\frac{h}{1+\eta}
\right)^{\frac{p}{2}+1}
\int_{B(0,r\sqrt\frac{1+\eta}{h})^c}|y|^2e^{-\cin{y}}\dr y
}{
\int_{B(0,r)}e^{-\cinh{v}}\dr v
+
\left( \frac{h}{1+\eta} \right)^\frac{p}{2}
\int_{B\left( 0,r\sqrt{\frac{1+\eta}{h}} \right)^c}
e^{-\cin{y}}\dr y
}   =  +\infty  .  \eqno (A.19)$$
Since the function $\bar\d$ is continuous, the last two formulas imply step 1. 

\noindent{\bf Step 2.} We refine step 1. Namely, we want to show that there are $\b>0$ and $h_0(r,\d)>0$ such that, if $(\l,\eta)$ is the couple of step 1 and $h\in(0,h_0(r,\d))$, then 
$$-1+h\log \frac{1}{h^{\b}}\le
\eta\le -1+\sqrt h   .   \eqno (A.20)$$
To show this, we let  
$$\a_-(h)=h
\log\frac{1}{h^{\b}}             ,\qquad 
\a_+(h)=\sqrt h  ,  $$
$$\eta_-(x)=-1+\a_-(h),\qquad \eta_+(x)=-1+\a_+(h)  .  $$
As in step 1, the assertion follows by the intermediate value theorem if we show that, for $h$ small, the function $\bar\d$ of step 1 satisfies
$$\bar\d(\eta_+(h))\le\d\le\bar\d(\eta_-(h))  .  \eqno (A.21)$$
Actually, we are going to show that that 
$\bar\d(\eta_+(h))\tends 0$ and $\bar\d(\eta_-(h))\tends +\infty$ for $h\tends 0$. 
To show this, we rewrite $\bar\d(\eta)$ as in (A.19). 
$$\bar\d(-1+\a_\pm(h))=$$
$$\2\left(
\frac{h}{\a_\pm(h)}
\right)^\frac{p}{2}       \cdot
\frac{1}{\a_\pm(h)}     \cdot
\frac{
\int_{B\left( 0,r\sqrt\frac{\a_\pm(h)}{h} \right)^c}
|y|^2e^{-\cin{y}}    \dr y
}{
\int_{B(0,r)}e^{-\cinh{v}}\dr v+
\left( \frac{h}{\a_\pm(h)} \right)^\frac{p}{2}
\int_{B\left( 0,r\sqrt\frac{\a_\pm(h)}{h} \right)^c}
e^{-\cin{y}}    \dr y
}    .   \eqno (A.22)$$Note that
$$\int_{
B\left(
0,r\sqrt\frac{\a_+(h)}{h}
\right)^c
}
|y|^2  e^{-\cin{y}}\dr y=
\int_{B\left(0,\frac{r}{h^\frac{1}{4}}\right)^c}
|y|^2e^{-\cin{y}}\dr y  \le
e^{
-\frac{D_4}{\sqrt h}
}  .  $$
By the last two formulas, the definition of $\a_+$ and the right hand side of (A.17) we get the inequality below, while (A.18) implies the limit.
$$\bar\d(-1+\a_+(h))\le
\2\cdot h^\frac{p}{4}\cdot\frac{1}{{\sqrt h}}\cdot
\frac{e^{-\frac{D_4}{\sqrt h}}}{(2\pi h)^\frac{p}{2}[1-a(h)]}  
\tends 0\txt{as}h\tends 0  .  $$
This yields the left half of (A.21). 

Now we show the other half of (A.21), i. e. that 
$\bar\d(\eta_-(h))\tends +\infty$ for $h\tends 0$. 
The first equality below comes from the definition of $\a_-(h)$, the second one from spherical coordinates; if $\e>0$ is given, the first and last inequalities hold if $h$ is small enough.
$$\int_{B\left( 0,r\sqrt\frac{\a_-(h)}{h} \right)^c}
e^{-\cin{y}}\dr y=
\int_{B\left( 0,r\sqrt{  \log\frac{1}{h^\b}  } \right)^c}
e^{-\cin{y}}\dr y=
C\int_{r\sqrt{  \log\frac{1}{h^\b}  }}^{+\infty}
\r^{p-1}e^{-\cin{\r}}\dr\r\le$$
$$C\int_{r\sqrt{  \log\frac{1}{h^\b}  }}^{+\infty}
e^{-\frac{\r^2}{2+\e}}\dr\r\le
C
\int_{r\sqrt{  \log\frac{1}{h^\b}  }}^{+\infty}
e^{-r\sqrt{  \log\frac{1}{h^\b}  }\cdot\frac{\r}{2+\e}}\dr\r=$$
$$-\left.
\frac{2+\e}{r\sqrt{\log\frac{1}{h^\b}}}\cdot
\exp\left\{
-\frac{r\sqrt{\log\frac{1}{h^\b}}}{2+\e}\r
\right\}     
\right\vert_{r\sqrt{\log\frac{1}{h^\b}}}^{+\infty}=
\frac{2+\e}{r\sqrt{\log\frac{1}{h^\b}}}
\exp\left\{
-\frac{r^2}{2+\e}\log\frac{1}{h^\b}
\right\}   \le
D_2 h^{\b\frac{r^2}{2+2\e}}   .  $$
Setting $D_6=\frac{r}{2+2\e}$, this yields the second inequality below; the second one is analogous; it is easy to check that 
$D_6\ge D_7$ and that $D_6-D_7\tends 0$ as $\e\tends 0$. 
$$h^{D_7\b}     \le
\int_{
B\left(
0,r\sqrt\frac{\a_-(h)}{h}
\right)^c
}
|y|^2e^{-\cin{y}}\dr y  ,\qquad
\int_{
B\left(
0,r\sqrt\frac{\a_-(h)}{h}
\right)^c
}
e^{-\cin{y}}\dr y\le h^{D_6\b}  .  \eqno (A.23)  $$
It may look a little strange that, in the formula above, $D_7\le D_6$ but recall that $\fun{}{x}{h^x}$ is monotone decreasing if 
$h\in(0,1)$ and that the first integral is larger than the second one. 

The last formula, (A.22), (A.17) and the definition of $\a_-$ yield 
$$\bar\d(-1+\a_-(h))\ge\2\cdot
\left(
\frac{1}{\log\frac{1}{h^\b}}
\right)^\frac{p}{2}\cdot
\frac{1}{ \log\frac{1}{h^\b} }\cdot
\frac{h^{D_7\b}}{
(2\pi h)^\frac{p}{2}+\left(
\frac{1}{\log\left(\frac{1}{h}\right)^\b}
\right)^\frac{p}{2}  h^{D_6\b}
}   .  \eqno (A.24)$$
We choose $\b>0$ so that
$$D_6\b>\frac{p}{2}  \txt{and}
D_7\b<\frac{p}{2}   .   \eqno (A.25)$$
This is possible since $D_7\le D_6$. By the first one of these formulas, the last term of the product on the right in (A.24) goes like $h^{D_7\b-\frac{p}{2}}$; this, the second one of (A.25) and (A.24) imply that 
$$\bar\d(-1+\a_-(h))\tends+\infty\txt{as}h\tends 0
\eqno (A.26)$$
ending the proof of (A.21). 

\noindent{\bf Step 3.} Let us call $\eta(\d,h)$ the solution of 
$\bar\d(\eta)=\d$. We assert that the thesis follows if we show that
$$(\l(\eta(\d,h))-1)-\eta(\d,h)\cdot\d-
\log\left(\frac{1}{2\pi h}\right)^\frac{p}{2}\ge
\frac{\d}{2}   \txt{for}\d\ge\d_0(r,h) .
\eqno (A.27)$$
Indeed, (A.15) implies the second equality below, while (A.14) implies the third one.
$$\int_{\R^p}A_h(\g_{\l(\eta(\d,h)),\eta(\d,h)},v)\dr v=
\int_{\R^p}\left[
\cinh{v}+\log\g_{\l(\eta(\d,h)),\eta(\d,h)}(v)
\right]   \g_{\l(\eta(\d,h)),\eta(\d,h)}(v)\dr v    =  $$
$$\int_{\R^p}\left[
\cinh{v}+(\l(\eta(\d,h))-1)-\cinh{v}-
\frac{\eta(\d,h)}{2h}|v|^2 1_{B(0,r)^c}(v)
\right]    \g_{\l(\eta(\d,h)),\eta(\d,h)}(v)\dr v   =$$
$$(\l(\eta(\d,h))-1)-\eta(\d,h)\cdot\d  .  $$
The thesis follows from this, (A.27) and the definition of $F_h$ in (A.13). 

\noindent{\bf Step 4.} We prove (A.27). We begin with the case in which $\eta$ satisfies (A.20); we have seen in the steps above that this covers all $\d$'s from $\bar\d(\eta_+(h))$ (which from now on we shall call $\d_0(r,h)$; we have seen that it tends to zero as $h\tends 0$) to a large $\d$ (i. e. $\bar\d(\eta_-(h))$) which tends to infinity as $h\tends 0$. 

$\d_0(r,h)$ (which tends to zero as $h\tends 0$) to a large $\d$ which tends to infinity as $h\tends 0$. 

The first equality below comes from our choice of $\l(\eta)$, the second one is (A.15), the third one is the change of variables $y=\sqrt\frac{1+\eta}{h}v$.
$$1=\int_{\R^p}\g_{\l(\eta),\eta}(v)\dr v=
e^{\l(\eta)-1}\left[
\int_{B(0,r)}e^{-\cinh{v}}\dr v+
\int_{B(0,r)^c}e^{-\frac{1+\eta}{2h}|v|^2}\dr v
\right]   =  $$
$$e^{\l(\eta)-1}\left[
\int_{B(0,r)}e^{-\cinh{v}}\dr v+
\left(
\frac{h}{1+\eta}
\right)^\frac{p}{2}
\int_{B(0,r\sqrt\frac{1+\eta}{h})^c}e^{-\cin{y}}\dr y
\right] .  $$
We define
$$b(\eta,h)\colon=
\left(\frac{1}{2\pi}\right)^\frac{p}{2}
\int_{B\left( 0,r\sqrt{\frac{1+\eta}{h}} \right)^c}
e^{-\cin{y}}\dr y      .  $$
We define $a(h)$ as in (A.17) and choose $\eta=\eta(\d,h)$; the last two formulas yield
$$e^{\l(\eta)-1}(2\pi h)^\frac{p}{2}\left[
1-a(h)+\left(\frac{1}{1+\eta}\right)^\frac{p}{2}
b(\eta,h)
\right]  =1    .  \eqno (A.28)$$
The first equality below is the definition of $b(\eta,h)$; the  inequality follows from (A.20) and (A.23); the limit comes from the first one of (A.25).
$$\left(
\frac{1}{1+\eta}
\right)^\frac{p}{2}b(\eta,h)=
\frac{1}{h^\frac{p}{2}}
\left(
\frac{h}{1+\eta}
\right)^\frac{p}{2}
\left(
\frac{1}{2\pi}
\right)^\frac{p}{2}
\int_{B\left( 0,r\sqrt\frac{1+\eta}{h} \right)^c}
e^{-\cin{y}}\dr y\le$$
$$\frac{1}{h^\frac{p}{2}}
\left(\frac{1}{\log\frac{1}{h^\b}}\right)^\frac{p}{2}
\left(\frac{1}{2\pi}\right)^\frac{p}{2}
h^{D_6\b}\tends 0   .  $$
From this and (A.28) we get 
$$\left\vert
e^{\l(\eta(\d,h))-1}\cdot(1-a(h))
(2\pi h)^\frac{p}{2}-1
\right\vert\tends 0\txt{as}h\tends 0 . $$
Taking logarithms in the last formula and using (A.18), we get that there is $\e(h)\tends 0$ as $h\tends 0$ such that
$$\left\vert
\l(\eta(\d,h))-1-\log\left(
\frac{1}{2\pi h}
\right)^\frac{p}{2}
\right\vert   \le\e(h)  .  $$
This implies the first inequality below; the second one follows from two facts: first,  $-\eta(\d,h)\ge\frac{3}{4}$; this is true for $h$ small by the right hand side of (A.20). The second one is 
$\e(h)\le\frac{\d^2}{4}$; since $\e(h)\tends 0$ for $h\tends 0$, this is true if $\d\le\d_0(r,h)$ with $\d_0(r,h)\tends 0$ as $h\tends 0$.  
$$(
\l(\eta(\d,h))-1
)    -
\eta(\d,h)\cdot\d   -
\log\left(
\frac{1}{2\pi h}
\right)^\frac{p}{2}\ge
-\eta(\d,h)\cdot\d-\e(h)\ge\frac{\d}{2}   .  $$
This is (A.26) when (A.20) holds; we saw in (A.26) that 
$\bar\d(\eta^-(h))\tends+\infty$ when $h\tends 0$. 

We only sketch the case in which $-1<\eta\le \eta_-(h)$, i. e. the case in which $\bar\d(\eta)$ is really large. The proof is divided into two cases: $-1+h\le\eta\le-1+h\log\frac{1}{h^\b}$ and 
$-1<\eta\le-1+h$. We sketch the second one. 

In this case we use again (A.19) and we see that
$$\bar\d(\eta)\simeq\left(
\frac{1}{1+\eta}
\right)  .  $$
On the other hand, the denominator of (A.19) is $e^{1-\l(\eta)}$; taking logarithms, we see that
$$1-\l(\eta)\simeq\log\left(
\frac{h}{1+\eta}
\right)^\frac{p}{2}  .  $$
Recalling that $-1<\eta\le \eta_-(h)$ and that $\eta^-(h)\tends-1$, we see from the last two formulas that
$$\l(\eta)-1-\eta\cdot\bar\d(\eta)\ge
\frac{\bar\d(\eta)}{2}  $$
as we wanted.

\fin




\vskip 2pc
\centerline{\bf Bibliography}


\noindent [1] S. Adams, N. Dirr, M. Peletier, J. Zimmer, From the large deviation principle to the Wasserstein gradient flow: a new micromacro passage, Comm. Math. Phys., {\bf 307}, 791-815, 2011.

\noindent [2] L. Ambrosio, N. Gigli, G. Savar\'e, Gradient Flows, Birkhaeuser, Basel, 2005.

\noindent [3] L. Ambrosio, N. Gigli, G. Savar\'e, Heat flow and calculus on metric measure spaces with Ricci curvature bounded below - the compact case. Analysis and numerics of PDE's, Springer-INdAM Ser., {\bf 4}, 63-115, Springer, Milan, 2013. 

\noindent[4] L. Ambrosio, G. Savar\'e, L. Zambotti, Existence and stability for Fokker-Planck equations with log-concave reference measures, probability Theory and Related Fields, {\bf 145}, 517-564, 2009.

\noindent [5] P. Bernard, B. Buffoni, Optimal mass transportation and Mather theory, J. Eur. Math. Soc., {\bf 9}, 85-121, 2007.

\noindent [6] U. Bessi, A time step approximation scheme for a viscous version of the Vlasov equation, Advances in Mathematics, {\bf 266}, 17-83, 2014. 

\noindent [7] V. Bogachev, G. Da Prato, M. R\"ockner, Uniqueness for solutions of Fokker-Planck equations on infinite-dimensional spaces, J. Evol. Equ., {\bf 10}, 487-509, 2010.

\noindent [8] G. Da Prato, F. Flandoli, M. R\"ockner, Uniqueness for continuity equations in Hilbert spaces with weakly differentiable drift, Partial Differ. Equ. Anal. Comput. {\bf 2}, 121-145, 2014.

\noindent [9] G. Da Prato, Introduction to stochastic analysis and Malliavin calculus, Pisa, 2007.




\noindent [10] L. De Pascale, M. S. Gelli, L. Granieri, Minimal measures, one-dimensional currents and the Monge-kantorovich problem, Calc. Var. and Partial Differential Equation, {\bf 27}, 363-388, 2006.

\noindent [11] I. Ekeland, R. Temam, Convex analysis and variational problems, Amsterdam, 1976.


\noindent [12] J. Feng, T. Nguyen, Hamilton-Jacobi equations in space of measures associated with a system of conservation laws, Journal de Math\'ematiques pures et Appliqu\'ees, {\bf 97}, 318-390, 2012. 






\noindent [13] D. Gomes, E. Valdinoci, Entropy penalization method for Hamilton-Jacobi equations, Advances in Mathematics, {\bf 215}, 94-152, 2007.



\noindent [14] C. L\'eonard, From the Schr\"odinger problem to the Monge-Kantorovich problem, Journal of Functional Analysis, {\bf 262}, 1879-1920, 2012.

\noindent [15] E. Nelson, Dynamical theories of Brownian motion, Princeton, 1967.

\noindent [16] C. Villani, Topics in optimal transpotation, Providence, R. I., 2003.

\end